\documentclass[A4paper,twoside,11pt,english]{article}
\usepackage{babel}
\usepackage{graphicx}
\usepackage{latexsym}
\usepackage{amsmath}
\usepackage{amsthm}
\usepackage{amsfonts}
\usepackage{mathrsfs}
\usepackage{geometry}
\usepackage[latin1]{inputenc}
\usepackage{amssymb}
\usepackage{makeidx}
\usepackage[table,xcdraw]{xcolor}
\usepackage{makeidx}
\usepackage{tocbibind}
\usepackage{tcolorbox}
\usepackage[framemethod=TikZ]{mdframed}
\usepackage{lipsum}
\usepackage[mathscr]{eucal}
\usepackage{array}
\usepackage{comment}
\usepackage{cite}
\usepackage{hyperref}
\usepackage{algorithm2e}
\usepackage{adjustbox}
\usepackage{caption}
\usepackage{subcaption}
\usepackage{comment}
\usepackage{multirow}
\usepackage{wrapfig}    % per testo attorno a immagini

\title{\textbf{Sampling Kantorovich operators for speckle noise reduction and gap filling with some applications to remote sensing}}
\author{Danilo Costarelli\thanks{Corresponding author. Paper submitted to arXiv (arXiv:2505.02422v2)} \quad - \quad Mariarosaria Natale\\
Department of Mathematics and Computer Science\\
University of Perugia\\
1, Via Vanvitelli, 06123 Perugia, Italy\\
{\tt danilo.costarelli@unipg.it} \& {\tt mariarosaria.natale@unipg.it}}
\date{}
\providecommand{\U}[1]{\protect\rule{.1in}{.1in}}
\theoremstyle{plain}
\newtheorem{theorem}{Theorem}

\newtheorem{corollary}[theorem]{Corollary}

\newtheorem{definition}[theorem]{Definition}

\newtheorem{lemma}[theorem]{Lemma}

\theoremstyle{definition}
\newtheorem{remark}[theorem]{Remark}

\newcommand{\R}{\mathbb{R}}
\newcommand{\Z}{\mathbb{Z}}

\newcommand{\be}{\begin{equation}}
\newcommand{\ee}{\end{equation}}

\newcommand{\supp}{\text{supp}\,}

\begin{document}
\maketitle              % typeset the header of the contribution
\begin{abstract}
In this paper, we investigate the application of multivariate sampling Kantorovich (SK) operators for image reconstruction, with a particular focus on gap filling and speckle noise reduction.
To understand the accuracy performances of the proposed algorithms, we first derive a quantitative estimate in $C(\R^n)$ for the error of approximation using the Euler-Maclaurin summation formula, under weak regularity conditions. We also establish a convergence result and a quantitative estimate with respect to the dissimilarity index measured by the continuous SSIM for functions in Lebesgue spaces. Additionally, we prove a multidimensional linear prediction result, which is used to design a new SK-based reconstruction algorithm to handle missing data, that we call LP-SK algorithm. To address speckle noise, we integrate SK operators into a newly proposed Down-Up scaling approach.
Numerical tests are presented on synthetic and real SAR images to validate the proposed methods. Performance is assessed using similarity metrics such as SSIM and PSNR, along with speckle-specific indexes. Comparative analysis with state-of-the-art techniques highlights the effectiveness of the proposed approaches.
\end{abstract}

\noindent {\small {\bf Keywords:} Sampling Kantorovich operators; Linear prediction; Gap filling; Missing information; Speckle noise reduction; SAR image; Quantitative estimates.

\section{Introduction}\label{s1}
Image restoration problems, such as gap filling and noise reduction, play a central role
in many imaging applications, including remote sensing, medical imaging, and signal reconstruction. Among these applications, the enhancement of remote sensing (RS) images has gained increasing interest due to their crucial role in Earth Observation (EO) applications. RS images, particularly those obtained from Synthetic Aperture Radar (SAR) and GNSS Reflectometry (GNSS-R) \cite{Zavo}, are often affected by degradations, that compromise the accurate interpretation and analysis of the data. This is especially critical when monitoring biogeophysical parameters related to Essential Climate Variables (ECVs), including for instance soil moisture and forest biomass \cite{Pierdicca2010, Moreau}. These variables play a significant role in understanding climate change, hydrological processes, and land cover dynamics, with implications for agriculture, climatology, surface temperature regulation, and environmental risk management. Two common types of degradation in RS images are a multiplicative noise (commonly referred to as \textit{speckle}) and missing data (\textit{gap areas}). Speckle noise is a granular disturbance that affects SAR images and other coherent imaging techniques, such as ultrasound and laser imaging. It results from the coherent nature of the acquisition of radar signals and significantly degrades image quality. Gap areas, on the other hand, refer to partially missing data caused by various factors, including cloud coverage in optical imagery, shadowed regions in SAR data, and sensor malfunctions such as the SLC-off problem, which occurred after the failure of the Scan Line Corrector (SLC) on Landsat 7 ETM+ in May 2003, resulting in wedge-shaped data gaps across the imagery \cite{SLCoff}. Geometric inconsistencies due to differences between satellite orbits and Earth rotation also contribute to these gaps. As a result, the development of effective denoising and gap filling techniques is essential to improve the usability of RS data.\\

Motivated by the above considerations, the present paper presents two novel methods for image restoration tasks, namely gap filling and speckle noise reduction, with particular motivation from remote sensing applications. The first method deals with gap filling and is based on the use of multivariate sampling Kantorovich (SK) operators, which have been studied in recent years \cite{CVbumi, GHR2011, OT2016, GS2025, BBSV2007, ACLR2017, AAK2022, CG2017} and have proven effective for image reconstruction and resolution enhancement \cite{CSV}. Here, we prove a linear prediction result specifically for the multivariate SK operators, allowing signal reconstruction using only ``past samples'' (the so-called {\em linear prediction}). Actually, in the proposed model, the variables represent spatial coordinates $(x,y)$, and then ``past pixels'' could be called ``previous pixels'', i.e., pixels $(i, j)$ taken in a neighborhood of $(\nu,\mu)$ with $i<\nu$ and $j<\mu$. This result serves as the theoretical foundation for a new algorithm, referred to as the LP-SK algorithm. A key advantage of this approach is its reliance solely on past pixel values, without requiring future or global image data. This may make the method suitable for applications where data becomes available over time, such as in data streams or real-time image detection, where future or complete image data may not yet be available. Other studies have explored multi-source or hybrid methods that combine data from different satellite missions \cite{multigap1, multigap2, LLJS2023}. Although effective, these approaches are based on the availability of auxiliary data, which, as we stress, may not always be accessible. In contrast, the proposed method uses only single-source information and exploits the spatial structure of the image, making it more broadly applicable. \\
The second goal of this paper is to address the reduction of speckle noise. We emphasize that speckle noise poses significant challenges for image interpretation and analysis \cite{despeckle1, SX1996, WL2025speckle}. In addition to the standard filtering approach, we propose an alternative method that integrates the SK algorithm implemented in \cite{CVbumi} and introduces two additional steps: downscaling and upscaling. For this reason, we refer to it as the Down-Up processing algorithm. Between these two steps, the method applies a filtering phase using some filters known from the literature. Quantitative testing demonstrates a general improvement in image quality. In particular, the proposed Down-Up processing algorithm significantly improves without-reference metrics (i.e., metrics that do not rely on ground-truth noiseless images), which are especially valuable in real-world applications where noise-free reference images are typically unavailable.\\

The structure of the paper is as follows. 

Section \ref{s2} provides a theoretical analysis of SK operators. Over the past forty years, sampling-type series have been a central topic in Approximation Theory, due to their deep connections with signal and image processing. In this section, we focus on the case of non-uniformly sampled data and present a new quantitative approximation result in $C(\R^n)$ (Theorem \ref{th_nuovo}), derived via the Euler-Maclaurin summation formula. We also prove the convergence and approximation order of the proposed algorithm (Algorithm \ref{SKalgorithm}) in the setting of functions belonging to $L^1(\R^n)\cap L^2(\R^n)$, with respect to the continuous dissimilarity index defined via the cSSIM, recently introduced by Marchetti \cite{marchetti2021}, which represents the continuous version of the discrete SSIM (\textit{structural similarity index measurement}). The theoretical results regarding the approximation order of SK operators will be essential for evaluating the accuracy performance of the proposed algorithms. Finally, we provide a qualitative estimate for functions in a suitable Lipschitz class (Corollary \ref{corollary_Lipschitz}), which is closely related to the decay of numerical errors associated to the SSIM, computed in the subsequent sections.

Section \ref{s3} is devoted to the extension of SK-based reconstruction to the case of missing data. We first establish a multidimensional linear prediction result (Theorem \ref{th_predizione}), which forms the basis for the implementation of the LP-SK algorithm (Algorithm \ref{predictionfunction}) designed to fill gaps in single-source satellite images. In terms of reconstruction quality, the proposed approach seems to be both efficient and accurate. A key advantage of the LP-SK algorithm is its reliance solely on past values, without requiring access to future or full image data. The numerical testing in this section confirms the potential of the method for filling gaps in realistic contexts.

Section \ref{s4} addresses the problem of speckle noise. After an overview of its nature and effects on images, we propose a Down-Up processing algorithm that combines classical interpolation techniques with the use of SK operators for scaling. This method aims to preserve structural details while reducing speckle noise. We evaluate the performance of the filtering methods both with- and without-reference indexes with the proposed Down-Up scheme. The results are consistent with previous findings \cite{CSV, CAR}, suggesting that SK approximation is particularly effective for upscaling, whereas bicubic interpolation performs better during downscaling. Consequently, the most effective strategy appears to be a bicubic downscaling step followed by despeckle filtering and SK upscaling algorithm. Our experiments prove a consistent improvement in without-reference quality metrics.

Finally, Section \ref{s5} outlines possible directions for future research.
\section{SK operators: theoretical analysis and applications to uniform samples}\label{s2}
We recall the definition and main properties of the sampling Kantorovich operators, which form the basis of the SK algorithm.\\

Let $\Pi^{n}=(t_{\underline{k}})_{\underline{k}\in\mathbb{Z}^{n}}$ be a sequence of vectors defined by $t_{\underline{k}}=(t_{k_1},\cdots,t_{k_n})$, where each $(t_{k_i})_{k_i\in\mathbb{Z}}$, $i=1,\cdots ,n$, is a sequence of real numbers with $-\infty<t_{k_i}<t_{k_{i}+1}<+\infty$, $\lim_{k_i\to \pm \infty}t_{k_i}=\pm \infty$, for every $i=1,\cdots,n$ and such that there exist $ \Delta, \delta >0$ for which $\delta\leq \Delta_{k_i}:=t_{k_{i}+1}-t_{k_i}\leq \Delta$, for every $i=1,\cdots,n$. Moreover, we denote by $R_{\underline{k}}^{w}:=\left [ \frac{t_{k_1}}{w},\frac{t_{k_1+1}}{w} \right ]\times \cdots \times \left [ \frac{t_{k_n}}{w},\frac{t_{k_n+1}}{w} \right ]$, $w>0$, whose Lebesgue measure is given by $A_{\underline{k}}/w^n$, where $A_{\underline{k}}:=\Delta_{k_1}\cdot\Delta_{k_2}\cdots \Delta_{k_n}$.\\

A function $X:\mathbb{R}^{n}\to\mathbb{R}$ will be called \textit{kernel} if it satisfies the following conditions
\begin{description}
\item ($\chi$1) $X\in L^1(\R^n)$ and it is bounded in a neighborhood of $\underline{0}\in\R^n$;
\item ($\chi$2) $\sum_{\underline{k}\in\Z^n}X(\underline{x}-t_{\underline{k}})=1$, for every $\underline{x}\in\R^n$;
\item ($\chi$3) for some $\beta>0$, $m_{\beta,\Pi^n}(X):=\displaystyle\sup_{\underline{x}\in\mathbb{R}^n}\sum_{\underline{k}\in\mathbb{Z}^n} |X(\underline{x}-t_{\underline{k}})| \left \|  \underline{x}-t_{\underline{k}}\right \| _{2}^{\beta}<+\infty$;
\end{description}
where $\left\|\cdot\right\|_{2}$ denotes the usual Euclidean norm on $\R^n$.
\begin{lemma}[Lemma 3.2 of \cite{CVbumi}]
Let $X$ be a kernel, then $m_{0,\Pi^n}(X):=\displaystyle\sup_{\underline{x}\in\mathbb{R}^n}\sum_{\underline{k}\in\mathbb{Z}^n} |X(\underline{x}-t_{\underline{k}})|<+\infty$. 
\end{lemma}

Many, but not all, multivariate kernels are defined as the $n$-fold (tensor) product of univariate kernel functions, such as the Fej\'{e}r kernel, the central B-splines, and others \cite{CVbumi, CCNV2022, CG2017}.
For instance, let $\chi\in L^{1}(\R)$ be a one-dimensional kernel (i.e., a kernel satisfying ($\chi$1)--($\chi$3) for $n=1$). We can define 
\begin{equation}\label{def_kernel_prod}
X(\underline{x}):=\prod_{i=1}^{n} \chi(x_i),\qquad \underline{x}\in\R^n.
\end{equation}
Since $\chi$ is integrable, it follows that $X\in L^{1}(\R^n)$, as shown by the following computation
\begin{equation*}
	\int_{\R^n}X(\underline{x})d	\underline{x}=\int_{\R^n}\prod_{i=1}^{n}\chi (x_i)dx_1\cdots dx_n=\prod_{i=1}^{n}\int_{\R}\chi(x_i)dx_i<+\infty,
\end{equation*}
and
\begin{equation*}
\sum_{\underline{k}\in\Z^n}X(\underline{x}-t_{\underline{k}})=\prod_{i=1}^{n}\sum_{k_i\in\Z}\chi(x_i-t_{k_i})=1,
\end{equation*}
for every $\underline{x}\in\R^n$. Furthermore, $m_{0,\Pi^n}(X)=\prod_{i=1}^{n}m_{0,\Pi}(\chi)<+\infty$, where the convergence of the involved series is uniform on the compact subsets of $\R^n$. 

\begin{definition}
The multivariate sampling Kantorovich operators for a given kernel $X$ are defined by
\begin{equation*}
(K_w^{X} f)(\underline{x}):=\sum_{\underline{k}\in\mathbb{Z}^n}X( w\underline{x}-t_{\underline{k}})\frac{w^n}{A_{\underline{k}}}\int_{R_{\underline{k}}^{w}}f(\underline{u})d\underline{u}, \qquad \underline{x}\in\mathbb{R}^n,	
\end{equation*}
where $f:\mathbb{R}^n\to\mathbb{R}$ is a locally integrable function such that the above series is convergent for every $\underline{x}\in\mathbb{R}^n$.
\end{definition}
The family of sampling Kantorovich operators provides one of the possible generalizations of the classical Whittaker-Kotel'nikov-Shannon (WKS) theorem, in which it has been proved that a signal can be reconstructed from its pointwise samples combined with the sinc function. A first well-known generalization of this scheme was introduced by Butzer and his coauthors (see, e.g., \cite{BFS1990}), leading to the family of generalized sampling operators. However, this approach remains sensitive to discontinuities, noise, and timing errors, since it relies on exact values of $f(\underline{k}/w)$. The sampling Kantorovich operators, introduced in \cite{BBSV2007}, distinguish themselves from the previous ones by replacing the sampled values of the signal $f$ with the integral mean on $R_{\underline{k}}^{w}$. This averaging provides a more robust approximation tool: it inherits
the regularization properties of the integral mean, which also helps in the reduction of time-jitter error. The approximation properties of sampling Kantorovich operators have been studied in several functional settings, including continuous functions, $L^p$-spaces and Orlicz spaces \cite{CVbumi, OT2016, CPV1,GS2025, ACLR2017}. Moreover, these operators act simultaneously as a low-pass filter, reducing noise, and as rescaling algorithm with better performances respect to other algorithms designed for the same purpose (see, e.g., \cite{CSV}).\\

In order to deal with the implementation of the sampling Kantorovich operators and its application to digital images, it can be useful to recall some key approximation results about their order of convergence.\\ 

For uniformly continuous and bounded functions $f\in C(\R^n)$, we refer to the well-known modulus of continuity, defined as $\omega(f,\delta):=\sup_{\left \| \underline{t} \right \|_{2}\leq \delta} |f(\cdot+\underline{t})-f(\cdot)|$, with $\delta>0$.
The estimate in the following theorem depends on whether condition ($\chi$3) holds with $\beta \geq 1$ or $0 < \beta < 1$. Certain kernels, such as the Fej\'{e}r kernel (see, e.g., \cite{CG2017}), do not have finite discrete absolute moments for $\beta \geq 1$ (more precisely, the first order discrete absolute moment is infinite), but satisfy condition ($\chi$3) for every $0 < \beta < 1$. Therefore,  we state the theorem by considering these two cases separately.
\begin{theorem}[Theorem 4.3 and Theorem 4.4 of \cite{CNV2025}]\label{th_stimeC}
Let $X$ be a kernel, and let $f\in C(\R^n)$. Then, we have the following inequalities, for sufficiently large $w>0$
\begin{enumerate}
\item if $X$ satisfies condition ($\chi$3) with $0<\beta< 1$:
\begin{equation*}
\left\|K^X_w f-f\right\|_{\infty}\leq C_1\omega\left(f, w^{-\beta}\right)+ 2^{\beta+1}\left\|f\right\|_{\infty}\;w^{-\beta}\;m_{\beta,\Pi^n}(X),
\end{equation*}
\item if $X$ satisfies condition ($\chi$3) with $\beta\geq 1$:
\begin{equation*}
\left\|K^X_w f-f\right\|_{\infty}\leq C_2 \omega\left(f,w^{-1}\right);
\end{equation*}

\end{enumerate}
where $C_1:=m_{0,\Pi^n}(X)+m_{\beta,\Pi^n}(X)+\,n^{\beta/2}\Delta^{\beta}m_{0,\Pi^n}(X)$, and $C_2:=m_{0,\Pi^n}(X)+\sqrt{n}\Delta m_{0,\Pi^n}(X)+m_{1,\Pi^n}(X)$.
\end{theorem}
It is well-known that, if ($\chi$3) is satisfied for $\beta\geq 1$, then it also holds for every values between $0$ and $\beta$.

Now, we focus on the case 1. of Theorem \ref{th_stimeC}, where $m_{\beta,\Pi^n}(X)<+\infty$, for every $0<\beta<1$, but $m_{1,\Pi^n}(X)=+\infty$ (exactly when we are not in the optimal order case), in order to derive the most possible accurate quantitative estimate for the aliasing error. We follow the approach used in the recent paper \cite{CC2025} in the univariate case, which relies on a special form of the Euler-Maclaurin summation formula.
\begin{theorem}[Formula 2.1.2 of \cite{EMS}]\label{th_EMS}
Let $f : [a, b] \to \mathbb{R}$ be a function whose derivative exists and is continuous on  $[a, b]$. Then  
\begin{equation*}
\sum_{a<n\leq b}f(n)=\int_{a}^{b} f(t)dt +\int_{a}^{b} \{t\} f'(t) dt - \{b\} f(b) + \{a\} f(a),
\end{equation*}
where the fractional part function is defined as $\{x\} := x - \lfloor x \rfloor \in [0, 1)$, for any real number $x$.
\end{theorem}
Now, we can state the following.
\begin{theorem}\label{th_nuovo}
Let $X$ be a kernel defined as in (\ref{def_kernel_prod}), such that 
\begin{equation}\label{ipotesi_chi}
\chi(x)=\mathcal{O}(|x|^{-2}), \text{ as } |x|\to+\infty, \text{ and with } m_{1,\Pi}(\chi)=+\infty.
\end{equation}
Then, for every $f\in C(\R^n)$, we have
\begin{equation*}
\left\|K^X_w f-f\right\|_{\infty}\leq K \left[\omega\left(f, \frac{\log w}{w}\right)+ \frac{\left\|f\right\|_{\infty}}{w}\right],
\end{equation*}
for every sufficiently large $w >0$, where $K >0$ is a suitable absolute constant depending only on the kernel $X$.
\end{theorem}
\begin{proof}
Since $\chi$ satisfies (\ref{ipotesi_chi}), there exists constants $\frac{x_0}{\sqrt{n}}>0$ and $C>0$ such that $|\chi(x)|\leq C|x|^{-2}$, for all $|x|>\frac{x_0}{\sqrt{n}}$. Let $w>x_0>\frac{x_0}{\sqrt{n}}$ and $\underline{x}\in\R^n$ be fixed. Using property ($\chi$2), we have
\begin{equation*}
\begin{split}
|(K_w^{X}f)(\underline{x})-f(\underline{x})|&\leq \sum_{\underline{k}\in\mathbb{Z}^n} |X( w\underline{x}-t_{\underline{k}})| \frac{w^n}{A_{\underline{k}}}\int_{R_{\underline{k}}^{w}}|f(\underline{u})-f(\underline{x})|d\underline{u}\\
&=\biggl\{\sum_{\left\|w\underline{x}-t_{\underline{k}}\right\|_{2}\leq w}+\sum_{\left\|w\underline{x}-t_{\underline{k}}\right\|_{2}> w}\biggr\}|X( w\underline{x}-t_{\underline{k}})| \frac{w^n}{A_{\underline{k}}}\int_{R_{\underline{k}}^{w}}|f(\underline{u})-f(\underline{x})|d\underline{u}\\
&=:I_1+I_2.
\end{split}
\end{equation*}
We estimate $I_1$. Exploiting the well-known inequality $\omega(f,\lambda\delta)\leq (1+\lambda)\omega(f,\delta)$, with the specific choice of $\lambda=\frac{w}{\log w}\left\|\underline{u}-\underline{x}\right\|_{2}$ and $\delta=\frac{\log w}{w}$, we get
\begin{equation*}
\begin{split}
I_{1}&\leq \sum_{\left\|w\underline{x}-t_{\underline{k}}\right\|_{2}\leq w} |X( w\underline{x}-t_{\underline{k}})| \frac{w^n}{A_{\underline{k}}}\int_{R_{\underline{k}}^{w}}\omega(f,\left\|\underline{u}-\underline{x}\right\|_{2}) d\underline{u}\\
&\leq \omega\left(f,\frac{\log w}{w}\right)\left[m_{0,\Pi^n}(X)+\frac{w}{\log w}\sum_{\left\|w\underline{x}-t_{\underline{k}}\right\|_{2}\leq w} |X( w\underline{x}-t_{\underline{k}})| \frac{w^n}{A_{\underline{k}}}\int_{R_{\underline{k}}^{w}}\left\|\underline{u}-\underline{x}\right\|_{2}d\underline{u}\right].
\end{split}
\end{equation*}
Now, for every $\underline{x},\underline{u}\in\R^n$, we may write
\begin{equation*}
\left\|\underline{u}-\underline{x}\right\|_{2}\leq \left\|\underline{u}-\frac{t_{\underline{k}}}{w}\right\|_{2}+\left\|\frac{t_{\underline{k}}}{w}-\underline{x}\right\|_{2}\leq \sqrt{n}\frac{\Delta}{w}+\frac{\left\|t_{\underline{k}}-w\underline{x}\right\|_{2}}{w},
\end{equation*}
for every $w>0$; therefore
\begin{equation*}
\begin{split}
I_1&\leq \omega\left(f,\frac{\log w}{w}\right)\left[m_{0,\Pi^n}(X)+\frac{\sqrt{n}\Delta}{\log w}m_{0,\Pi^n}(X)+\frac{1}{\log w}\sum_{\left\|w\underline{x}-t_{\underline{k}}\right\|_{2}\leq w} |X( w\underline{x}-t_{\underline{k}})|\left\|t_{\underline{k}}-w\underline{x}\right\|_{2}\right]\\
&\leq \omega\left(f,\frac{\log w}{w}\right)\left[m_{0,\Pi^n}(X)\left(1+\frac{\sqrt{n}\Delta}{\log w}+\frac{x_0}{\log w}\right)+\frac{1}{\log w}\sum_{x_0<\left\|w\underline{x}-t_{\underline{k}}\right\|_{2}\leq w} |X( w\underline{x}-t_{\underline{k}})|\left\|t_{\underline{k}}-w\underline{x}\right\|_{2}\right].
\end{split}
\end{equation*}
Recalling that in $\R^n$ all norms are equivalent, we specifically have the bound $$\left\|\underline{x}\right\|_{max}\leq \left\|\underline{x}\right\|_{2}\leq \sqrt{n}\left\|\underline{x}\right\|_{max},$$ for any $\underline{x}\in\R^n$, where $\left\|\underline{x}\right\|_{max}:=\max_{i=1,\cdots,n} |x_i|$. Using this, we derive the following inequality
\begin{equation*}
\begin{split}
&\sum_{x_0<\left\|w\underline{x}-t_{\underline{k}}\right\|_{2}\leq w} |X( w\underline{x}-t_{\underline{k}})|\left\|t_{\underline{k}}-w\underline{x}\right\|_{2}\leq \sqrt{n} \sum_{\frac{x_0}{\sqrt{n}}<\left\|w\underline{x}-t_{\underline{k}}\right\|_{max}\leq w} |X( w\underline{x}-t_{\underline{k}})|\left\|t_{\underline{k}}-w\underline{x}\right\|_{max}\\
&\leq \sqrt{n} \sum_{i=1}^{n}
\left\{\sum_{\frac{x_0}{\sqrt{n}}<|wx_i-t_{k_i}|\leq w}|\chi(wx_i-t_{k_i})||t_{k_i}-wx_i|\left[\sum_{\underline{k}_{[i]}\in\Z^{n-1}}X^{[i]}\left(w\underline{x}_{[i]}-t_{\underline{k}_{[i]}}\right)\right] \right\}\\
&=\sqrt{n} \sum_{i=1}^{n}\left\{\sum_{\frac{x_0}{\sqrt{n}}<|wx_i-t_{k_i}|\leq w}|\chi(wx_i-t_{k_i})||t_{k_i}-wx_i|\right\}=:\sqrt{n}\sum_{i=1}^{n}I_{1,i},
\end{split}
\end{equation*}
where $X^{[i]}(w\underline{x}_{[i]}-t_{\underline{k}_{[i]}}) := \chi(wx_1-t_{k_1}) \cdot \, 
        ... \, \cdot \chi(wx_{i-1}-t_{k_{i-1}}) \cdot 
    \chi(wx_{i+1}-t_{k_{i+1}}) \cdot \, ... \, \cdot \chi(wx_n-t_{k_n})$, with $\underline{x}_{[i]} := (x_1, ..., x_{i-1}, x_{i+1}, ..., x_n)$ $\in $$\R^{n-1}$, $t_{\underline{k}_{[i]}} := (t_{k_1}, ..., t_{k_{i-1}}, t_{k_{i+1}}, ..., t_{k_n}) \in \Z^{n-1}$, for every $i= 1,...,n$.\\
Now, we focus on $I_{1,i}$ and, using (\ref{ipotesi_chi}), we get
\begin{equation*}
\begin{split}
I_{1,i}&\leq C  \sum_{\frac{x_0}{\sqrt{n}}<|wx_i-t_{k_i}|\leq w}|t_{k_i}-wx_i|^{-1}\\
&= C \sum_{wx_i-w\leq t_{k_i}<wx_i-\frac{x_0}{\sqrt{n}}}(wx_i-t_{k_i})^{-1}+C \sum_{wx_i+\frac{x_0}{\sqrt{n}}< t_{k_i}\leq wx_i+w}(t_{k_i}-wx_i)^{-1}\\
&=C \sum_{\frac{x_0}{\sqrt{n}} -wx_i< \widetilde t_{k_i}\leq w-wx_i}(wx_i + \widetilde t_{k_i})^{-1}+C \sum_{wx_i+\frac{x_0}{\sqrt{n}}< t_{k_i}\leq wx_i+w}(t_{k_i}-wx_i)^{-1}.
\end{split}
\end{equation*}
We can apply Theorem \ref{th_EMS} to the functions $F_1(t):=1/(wx_i+t)$ and $F_2(t):=1/(t-wx_i)$, which are differentiable with a continuous derivative on the intervals $[\frac{x_0}{\sqrt{n}} -wx_i, w-wx_i]$ and $[wx_i+\frac{x_0}{\sqrt{n}}, wx_i+w]$, respectively. As a result, we can write $I_{1,i}\leq C(C_1+ C_2)\log w$, for sufficiently large $w>0$, where the constants $C_1, C_2$ are absolute.\\
Then we can conclude that
\begin{equation*}
\begin{split}
I_{1}&\leq \omega\left(f,\frac{\log w}{w}\right)\left[m_{0,\Pi^n}(X)\left(1+\frac{\sqrt{n}\Delta}{\log w}+\frac{x_0}{\log w}\right)+ n\sqrt{n}C(C_1+C_2)\right]\leq K_1 \omega\left(f,\frac{\log w}{w}\right),
\end{split}
\end{equation*}
for $w>0$ sufficiently large, and a suitable absolute constant $K_1$ (independent of $w$). Finally, for what concerns $I_2$ we have
\begin{equation*}
\begin{split}
I_{2}&\leq 2\left\|f\right\|_{\infty} \sum_{\left\|w\underline{x}-t_{\underline{k}}\right\|_{2}> w}|X( w\underline{x}-t_{\underline{k}})| \leq 2\left\|f\right\|_{\infty} \sum_{\left\|w\underline{x}-t_{\underline{k}}\right\|_{max}> {\frac{w}{\sqrt{n}}}}|X( w\underline{x}-t_{\underline{k}})|\\
&\leq 2\left\|f\right\|_{\infty}\sum_{i=1}^{n}\left\{\sum_{|wx_i-t_{k_i}|>\textcolor{black}{\frac{w}{\sqrt{n}}}}|\chi (wx_i-t_{k_i})|\right\},
\end{split}
\end{equation*}
and by using again the assumption (\ref{ipotesi_chi}), we can write
\begin{equation*}
\begin{split}
I_{2}& \leq 2C\left\|f\right\|_{\infty}\sum_{i=1}^{n}\left\{\sum_{|wx_i-t_{k_i}|> {\frac{w}{\sqrt{n}}}}|wx_i-t_{k_i}|^{-2}\right\}\\
&=2C\left\|f\right\|_{\infty}\sum_{i=1}^{n}\left\{\sum_{t_{k_i}<wx_i-{\frac{w}{\sqrt{n}}}} (t_{k_i}-wx_i)^{-2}+\sum_{t_{k_i}>wx_i+{\frac{w}{\sqrt{n}}}} (t_{k_i}-wx_i)^{-2}\right\}\\
&\leq 8n\textcolor{black}{\sqrt{n}}C\left\|f\right\|_{\infty}\frac{1}{w}:=K_2\left\|f\right\|_{\infty}\frac{1}{w}.
\end{split}
\end{equation*}
The last inequality follows from the integral estimates
\begin{equation*}
\sum_{t_{k_i}<wx_i-\textcolor{black}{\frac{w}{\sqrt{n}}}} (t_{k_i}-wx_i)^{-2}\leq \int_{-\infty}^{wx_{i}-\textcolor{black}{\frac{w}{\sqrt{n}}+\Delta}} (t-wx_i)^{-2}dt=\frac{1}{\textcolor{black}{\frac{w}{\sqrt{n}}-\Delta}}\leq \frac{2\textcolor{black}{\sqrt{n}}}{w},
\end{equation*}
and
\begin{equation*}
\sum_{t_{k_i}>wx_i+\textcolor{black}{\frac{w}{\sqrt{n}}}} (t_{k_i}-wx_i)^{-2}\leq \int_{wx_i+\textcolor{black}{\frac{w}{\sqrt{n}}-\Delta}}^{+\infty} (t-wx_i)^{-2}dt=\frac{1}{\textcolor{black}{\frac{w}{\sqrt{n}}-\Delta}}\leq \frac{2\textcolor{black}{\sqrt{n}}}{w},
\end{equation*}
for $w>0$ sufficiently large. Thus, the thesis follows by taking $K= \max\{K_1, K_2\}$.
\end{proof}

In order to deal with not necessarily continuous functions (such as digital images), we consider $f\in L^{p}(\mathbb{R}^{n})$, $1\leq p<+\infty$. In this context, we can state the following theorem, after recalling the definition of the first-order $L^{p}$-modulus of smoothness of $f$, which is given by
\begin{equation*}
\omega(f, \delta)_{p}:=\sup_{\left \| \underline{t} \right \|_{2}\leq \delta} \left \| f(\cdot+\underline{t})-f(\cdot) \right \|_{L^p(\R^n)} =\sup_{\left \| \underline{t} \right \|_{2}\leq \delta} \left(\int_{\mathbb{R}^n}|f(\underline{s}+\underline{t})-f(\underline{s})|^{p}d\underline{s}\right)^{1/p}, \quad  \delta>0.
\end{equation*}
\begin{theorem}[Theorem 4.1 of \cite{CCNV2022}]\label{stimaLP}
Suppose that $M^{p}(X):=\displaystyle\int_{\mathbb{R}^n}|X(\underline{u})|\left \| \underline{u} \right \|_{2}^{p}d\underline{u}<+\infty$,
for some $1\leq p<+\infty$. Then, for every $f \in L^p (\mathbb{R}^n )$, the following quantitative estimate holds
\begin{equation*}
\left \| K^X_w f-f \right \|_{L^p(\R^n)}\leq C_{3} \omega\left ( f,\frac{1}{w} \right )_p ,
\end{equation*}
where $C_{3}:=\delta ^{-\frac{n}{p}}\left ( m_{0,\Pi^n}(X) \right )^{\frac{p-1}{p}}\left\{2^{\frac{p-1}{p}}m_{0,\Pi^n}(\tau)^{1/p}\left [ \left \| X\right \|_1 + M^{p} (X) \right ]^{\frac{1}{p}} + \left(m_{0, \Pi^n}(X)\right)^{\frac{1}{p}} \Delta^{\frac{n}{p}}(1+\sqrt{n}\Delta)\right\}$, for sufficiently large $w>0$, with $\tau$ denoting the characteristic function of the set $[0,\delta]^n$. 
\end{theorem}

Finally, several additional results have been obtained for the family of sampling Kantorovich operators. In recent years, these operators have been widely studied, leading to various generalizations, including nonlinear versions \cite{CNpreprint, CNV2025}.\\

Now, our aim is to establish a convergence result with respect to \textit{continuous structural similarity index} (cSSIM). This index was recently introduced in \cite{marchetti2021} (see also \cite{CNP1}) and represents the continuous version of the so-called (discrete) SSIM \cite{UIQIindex, SSIM1, SSIM2}, widely used in image analysis. To proceed, we recall its definition.\\

Let $\Omega\subset \R^n$ be a bounded domain and $\nu$ be a probability measure on $\Omega$. Let $f,g\in L^2(\Omega)$ be non-negative $\nu$-almost everywhere, the cSSIM is defined as follows
\begin{equation*}\label{def_cSSIM}
\text{cSSIM}(f,g) := \frac{2\mu_f\mu_g + c_1}{\mu_f^2 + \mu_g^2 + c_1}\cdot\frac{2\sigma_{fg} + c_2}{\sigma_f^2 + \sigma_g^2 + c_2},
\end{equation*}
with the constants $c_1, c_2 > 0$ are small stabilizing factors that prevent division by zero, where $\mu_f :=\left\|f\right\|_1$, $\sigma_{fg} :=\int_{\Omega}(f-\mu_f)(f-\mu_g)d\nu$, and $\sigma_f^2 := \sigma_{ff}$.\\
Although this definition assumes that $\Omega\subset \R^n$ is bounded, ensuring $L^2(\Omega)\subset L^1(\Omega)$ (in order to be well posed), it naturally extends to the entire space $\R^n$ by considering functions in $L^1(\R^n)\cap L^2(\R^n)$. In this more general setting, we can state the following result that can be immediately deduced as in Theorem 3.3 of \cite{marchetti_santin2022}.

\begin{theorem}\label{th_marchettisantin}
Let $f,g$ be two non-negative real-valued functions in $L^1(\R^n)\cap L^2(\R^n)$ and $c_1,c_2>0$. Then, the following bound holds
{\rm
\begin{equation*}
|1-\text{cSSIM}(f,g)|\leq c_f \left\|f-g\right\|_{L^2(\R^n)}^2,
\end{equation*}}
where $c_f:=\frac{4}{\sigma_f^2+c_2}+\frac{1}{\mu_f^2+c_1}$.
\end{theorem}
As a byproduct of Theorem \ref{stimaLP} and Theorem \ref{th_marchettisantin}, for any $f\in L^1(\R^n)\cap L^2(\R^n)$ we obtain the following estimate
\begin{equation*}
|1-\text{cSSIM}(f,K_w^X f)|\leq K \left(\omega \left ( f,w^{-1} \right )_2\right)^2,
\end{equation*}
as $w \to +\infty$, for a suitable $K>0$. If we further assume that $f$ belongs to a suitable Lipschitz class, we can immediately deduce the qualitative order of convergence. We recall that the Lipschitz classes in $L^p$-spaces for $0 < \alpha\leq 1$ are defined as $Lip (\alpha, p):=\{f\in L^{p}(\R^n) \;:\; \omega(f,\delta)_p=\mathcal{O}(\delta^{\alpha}), \;\mathrm{as} \;\, \delta\to 0 \}$.
\begin{corollary}\label{corollary_Lipschitz}
Under the assumptions of Theorem \ref{stimaLP}, and for any $f\in Lip(\alpha,2)\cap L^1(\R^n)$, with $0<\alpha\leq 1$, there exists a constant $\widetilde{K}>0$ such that
{\rm
\begin{equation*}
|1-\text{cSSIM}(f,K_w^X f)|\leq \widetilde{K} \frac{1}{w^{2\alpha}},
\end{equation*}}
for sufficiently large $w>0$.
\end{corollary}
\begin{proof}
Let $f$ be as in the statement. By Theorem \ref{stimaLP} with $p=2$, we obtain
$$
\| K_w^X f - f \|_{L^2(\R^n)}^2\, \leq C_3^2 \omega\left(f, \frac{1}{w}\right)^2_2 \leq C_3^2 K \frac{1}{w^{2\alpha}},
$$
for some constant $K>0$, and sufficiently large $w>0$. Then, applying Theorem \ref{th_marchettisantin}, the thesis follows by setting $\widetilde{K}=c_f\,C^2_3 K$.
\end{proof}

The family of sampling Kantorovich operators is revealed to be a suitable mathematical tool for the reconstruction and enhancement of images (referred to in what follows as SK algorithm), as it acts both as a low-pass filter and as a rescaling algorithm, increasing the resolution of images \cite{CSV}. For completeness, we provide the pseudo-code for the SK algorithm below (see Algorithm \ref{SKalgorithm}). Here, we consider the sequence $t_{\underline{k}} = \underline{k}$, $\underline{k} \in \mathbb{Z}^n$, and so we work with an equally spaced grid of nodes, where $\delta = \Delta = \Delta_{k_i} = 1$ and $A_{\underline{k}} = 1$ for every $\underline{k}\in \mathbb{Z}^n$.\\

Now, in order to apply the bi-dimensional  sampling Kantorovich operators to image reconstruction, we first recall how a discrete image can be identified with a piecewise constant function. Let $A=(a_{ij})_{i=1,\ldots,n;\,j=1,\ldots,m}$ denote a discrete grayscale image.
We associate with $A$ the step function, $I_A \in L^p(\mathbb{R}^2)$, $1 \leq p < +\infty $, with compact support, defined by
\begin{equation}\label{funz_imm}
I_A(x, y) := \sum_{i=1}^{n} \sum_{j=1}^{m} a_{ij} \cdot \mathbf{1}_{ij}(x, y), \quad (x, y) \in \mathbb{R}^2,    
\end{equation}
where
\begin{equation*}
\mathbf{1}_{ij}(x, y) := \begin{cases} 
1, & \text{if } (x, y) \in (i-1, i] \times (j-1, j], \\
0, & \text{otherwise},
\end{cases}   
\end{equation*}
is the characteristic function of the set $(i-1, i] \times (j-1, j]$, with $i=1,\cdots,n$ and $j = 1, \cdots, m $. 
In what follows, when no ambiguity arises, we refer to both $A$ and its continuous representation $I_A$ as an image.\\
Now, we can apply the bi-dimensional sampling Kantorovich operator to the function $I_A$, by selecting a suitable kernel function $X$. This allows us to reconstruct the continuous image associated with $A$ for any given sampling rate $w$ (\textit{sampling rate}). In particular, the parameter $w$ controls the density of the sampling grid: larger values of $w$ correspond to a finer sampling and therefore to a more accurate reconstruction of the original image structure. After applying the SK operator, the resulting function is evaluated on a discrete grid in order to obtain the final reconstructed digital image. This can be interpreted as a two-step process: first discretizations of the spatial domain (\textit{sampling}), and then discretizations of the range of intensity values (\textit{quantization}).

\vspace{0.3cm}
\RestyleAlgo{ruled}
\begin{algorithm}[H]\label{SKalgorithm}
\SetAlgoLined
\KwData{Original image $A$ ($n\times m$ pixel resolution); parameter $w>0$ of the operator $K_w^X$; integer scaling factor $r>0$.}
\KwResult{The reconstructed image of resolution $(n\cdot r)\times (m\cdot r)$.}
- Choice and definition of the kernel function $X$\;
- Computation of the size of the reconstructed image using the scaling factor: $(n\cdot r)\times (m\cdot r)$\;
- Modeling of $A$ as a function $I_A$ of the form (\ref{funz_imm})\;
- Definition of a sampling grid $G_r$, containing $(n\cdot r)\times (m\cdot r)$ uniformly spaced nodes over the square domain $[0,n]\times [0,m]$\;
- Computation of matrices of the sample values of $I_A$ by means of the Kronecker matrix product\;
- Definition of the vectors containing the arguments of $X$\;
\For{$i=1,\dots,n$ and $j=1,\dots,m$}{
sum over ${\underline{k}}$ of all non-zero or not negligible terms of the form $X(w\underline{x}_{ij}-\underline{k})\left[w^2\int_{R^{w}_{\underline{k}}}I_A(\underline{u})d\underline{u}\right]$, at the points $\underline{x}_{ij}$ belonging to the reconstruction grid $G_r$\;
}
\caption{Pseudocode of the SK algorithm}
\end{algorithm}
\vspace{0.3cm}

\begin{remark}
The parameter $r$ in Algorithm \ref{SKalgorithm} acts as a scaling factor for the resolution of the reconstructed image. The number of nodes in $G_r$ depends by the scaling factor $r$, and it is equal to $(n\cdot r) \times (m\cdot r)$ which corresponds to the final size of the rescaled image. In particular, if $0<r<1$ the algorithm performs a \textit{downscaling}, e.g., with $r=1/2$, the output image has size $\frac{n}{2}\times\frac{m}{2}$. If $r=1$, the reconstructed image preserves the original resolution. If $r>1$, the algorithm performs an \textit{upscaling}, increasing the resolution of the image by a factor of $r$.
\end{remark}

The performance of the SK algorithm in image rescaling, in terms of PSNR and CPU time, in comparison with some interpolation and quasi-interpolation methods, has been evaluated in \cite{CSV}.\\

A useful class of kernels for Algorithm \ref{SKalgorithm} is the \textit{Jackson type kernels} of order $s\in\mathbb{N}$, defined in the univariate case by
\begin{equation*}
J_{s}(x):=c_{s}\, \text{sinc}^{2s}\left(\frac{x}{2s\pi}\right), \;\;\; x\in\mathbb{R},	
\end{equation*}
where $c_{s}$ is a non-zero normalization coefficient, given by
\begin{equation*}
c_{s}:=\left[\int_{\mathbb{R}}\text{sinc}^{2s}\left(\frac{u}{2s\pi}\right) du\right]^{-1}.	
\end{equation*}
The multivariate Jackson-type kernel is given by the $n$-fold product of the corresponding univariate function, $\mathcal{J}_{s}^{n}(\underline{x})=\prod_{i=1}^{n} J_{s}(x_i)$, $\underline{x}\in\mathbb{R}^{n}$. For the numerical tests given in Section \ref{ss4_3}, we consider the bivariate
Jackson-type kernel with $s=12$. In fact, as stated in \cite{CSV}, in the case of rescaling images with double dimensions, it is sufficient to choose $w= 15$ when $s = 12$ to achieve ``good'' reconstructions.
\section{Implementation of the SK algorithm for missing data}\label{s3}
In this section, we propose a method to fill gaps in scattered data, a common issue in remote sensing techniques such as GNSS-R, where sampling is irregular. The method we propose uses SK operators, which were discussed in the previous section. This approach aims to improve the analysis of biogeophysical variables.\\

The theoretical foundation of the proposed algorithm (see Algorithm \ref{predictionfunction}, that in what follows will be referred to as the LP-SK algorithm) is based on the linear prediction of signals, using sampling Kantorovich operators, by sample values taken only from the ``past''. In general, Kantorovich sampling series approximate a signal at a given time by using sample values from both the ``past'' and the ``future'' (as happens, e.g., in the WKS--sampling theorem). However, this approach is not practical in real-world applications, where only past values of the signal are available. Here, we propose a multidimensional version of this result using a product-type kernel.

\begin{theorem}\label{th_predizione}
Let $X$ be a kernel defined as in (\ref{def_kernel_prod}), such that $\supp \chi\subset[\Delta,+\infty)$. Then, for every signal $f:\R^n\to\R$ for which the operators $(K_w^{X} f)_{w>0}$ are well-defined, and for any fixed $\underline{x}\in\mathbb{R}^n$, we have 
\begin{equation}\label{eq_prediction}
(K_w^{X} f)(\underline{x})=\sum_{\underline{k}\in\Z^n \atop t_{k_i}/w < x_i-\Delta/w}X( w\underline{x}-t_{\underline{k}})\,\frac{w^n}{A_{\underline{k}}} \int_{R_{\underline{k}}^{w}}f(\underline{u})d\underline{u},
\end{equation}
for every $w>0$.
\end{theorem}
\begin{proof}
Since $\supp \chi\subset[\Delta,+\infty)$, we have $\chi( wx_i-t_{k_i})=0$ for every $k_i\in\Z$ such that $t_{k_i}/w \geq x_i- \Delta/w$, $i=1,\cdots,n$. Hence, all integrals $\frac{w^n}{A_{\underline{k}}} \int_{R_{\underline{k}}^{w}}f(\underline{u})$ represent mean values of $f$ computed considering the values of $f$ on sets such that, for each dimension, $\frac{t_{k_i}}{w}, \frac{t_{k_i+1}}{w}<x_i$, $i=1,\cdots,n$.
\end{proof}

A key aspect of Algorithm \ref{predictionfunction} is the choice of the kernel function, which plays a fundamental role in the reconstruction process. For this purpose, we consider the well-known central B-spline (univariate) of order $s\in\mathbb{N}$, defined as
\begin{equation*}
B_{s}(x):=\frac{1}{(s-1)!}\sum_{j=0}^{s}(-1)^{j}\binom{s}{j}\left (\frac{s}{2} +x-j \right )_{+}^{s-1}\,,\quad x\in\R,
\end{equation*}
where $x_{+}:=\max\{ x,0 \}$ is the positive part of $x$. Its Fourier transform is $\widehat{B_{s}}(v)=\text{sinc}^{s}\left(\frac{v}{2\pi}\right)$, $v\in\mathbb{R}$, and it can be used to check that $B_s$ satisfies the partition of unity property, i.e., $\sum_{k\in\mathbb{Z}}B_{s}(u-k)=1$ (see \cite{CVbumi}). Moreover, $B_s$ is bounded on $\mathbb{R}$, with compact support on $[-s/2, s/2]$, ensuring their integrability in $L^{1}(\mathbb{R})$ with $\left\|B_s\right\|_{1}=1$. The bivariate central B-spline of order $s$ (that will be used in Algorithm \ref{predictionfunction}) is given by $\mathcal{B}_{s}^{2}(\underline{x}):=\prod_{i=1}^{2}B_s (x_i)$, $\underline{x}\in\mathbb{R}^{n}$.\\

In order to satisfy the assumptions of Theorem \ref{th_predizione}, the kernel $B_s$ must be shifted to the right so that its support is contained in $[\Delta, +\infty)$. In particular, in Algorithm \ref{predictionfunction}, we perform a translation by $\frac{s+2\Delta}{2}$. By the properties of the Fourier transform, condition ($\chi$2) is still satisfied.

\vspace{0.3cm}
\begin{algorithm}[H]\label{predictionfunction}
\SetAlgoLined
\KwData{Original image $A$ ($n\times m$); parameter $w>0$ of the operator $K_w^X$; coordinates $(\nu,\mu)$ of the missing pixel; order $s$ of B-spline kernel function.}
\KwResult{A value between 0 and 255, representing the predicted luminance of the missing pixel.}
- Modeling of $A$ as a function $I_A$ of the form (\ref{funz_imm})\;
- Computation of the sub-matrix of dimension $n_1\times m_1$ of the sample values of $A$ (each missing pixel is reconstructed using only past data in the sense specified in Remark \ref{remark_past})\;
- Definition of the mask of weights of dimension $n_1\times m_1$ determined by the shifted kernel function $\mathcal{B}_{s}^{2}$\;
\For{$i=1,\dots,n_1$ and $j=1,\dots,m_1$}{
sum over ${\underline{k}}$ of all non-zero of the form $\mathcal{B}_{s}^{2}(w\underline{x}-\underline{k})\left[w^2\int_{R^{w}_{\underline{k}}}I_A(\underline{u})d\underline{u}\right]$, for $\underline{x}=(\nu,\mu)$\;
}
\caption{Pseudocode of the LP-SK algorithm}
\end{algorithm}
\vspace{0.3cm}

\begin{remark}\label{remark_past}
Actually, in the proposed LP-SK algorithm, the variables represent spatial coordinates and do not involve any temporal sequences of data. Consequently, what we refer to as ``past pixels'' can more appropriately be called ``previous pixels'' (see Figure \ref{past_pixels}). For a missing pixel at coordinates $(\nu,\mu)$, we define as past pixels the samples available within a mask of size $n_1\times m_1$ located in the region $\{(i,j)\;:\; i<\nu,j<\mu\}$. This ensures that the reconstruction of each pixel depends only on previously scanned spatial data, in accordance with the assumptions of Theorem \ref{th_predizione}.  In the numerical experiments presented in this work, however, the notion of previous pixels is mainly used to show the local behaviour of the LP-SK algorithm which avoids the use of ``future samples''.
\end{remark}
\begin{figure}[htbp]
  \centering
\includegraphics[width=0.38\textwidth]{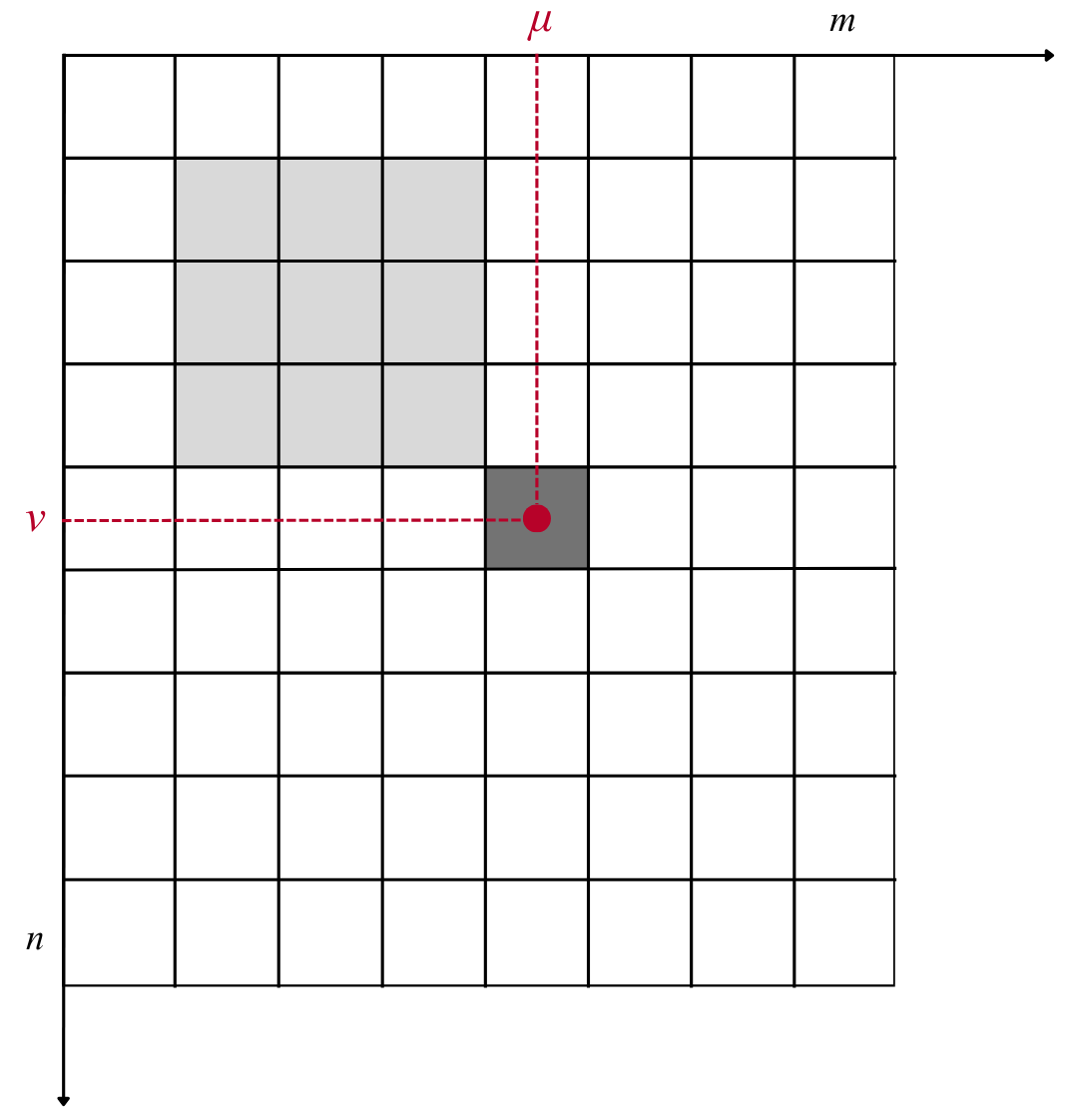}
  \caption{The dark gray cell indicates the missing pixel located at coordinates $(\nu,\mu)$. The light gray region represents the set of ``past pixels'' employed in the image reconstruction.}
\label{past_pixels}
\end{figure}

\subsection{Numerical testing}\label{ss3_1}
As mentioned before, gap filling in satellite images \cite{missing, aidareview} addresses the issue of missing data caused by several reasons, e.g. cloud coverage for optical imagery, shadowed area for SAR data sets, or instrumentation errors, e.g. SLC-off problem and line striping. Gap areas can have different sizes, dimensions, and locations. As a consequence, the usability of the data is limited. To address this problem, missing values are often excluded from the analysis or replaced by predictions from a variety of prediction methods (also called \textit{gap filling methods} or \textit{imputation methods}). These methods are typically categorized into three classes: \textit{single-source}, \textit{multi-source}, and \textit{hybrid methods}.

Multi-source methods \cite{multigap1, multigap2} rely on additional images, either from different dates or sensors, and generally offer more accurate reconstructions due to the increased amount of available information. In recent years, deep learning approaches \cite{deep, deep2, deep3} have gained significant traction in machine learning and are increasingly applied in remote sensing, because they can handle large datasets and learn complex patterns. Single-source methods, on the other hand, are based solely on information within the same image. These methods are particularly useful when no auxiliary data are available. Typically, single-source techniques \cite{single2, single1} are based on pixel similarity or dissimilarity rules, reconstructing missing areas by referencing non-gapped regions of the image. Common interpolation methods such as natural neighbor method \cite{NNN}, bilinear and cubic interpolation, Kriging interpolation method \cite{Kringing}, and inpainting algorithms \cite{inpainting, singlegap, adv1} are also used. Hybrid method is combination of single source and multi-source methods.\\

Within the single-source category, we propose a new method, which is based on SK linear prediction, referred to as the LP-SK algorithm (Algorithm \ref{predictionfunction}). \\
To evaluate its performance, we compare it with \textsc{Matlab\textsuperscript{\textcopyright}} inbuilt command function {\tt{imfill(A,conn)}}, which performs hole filling on grayscale images. In this context, a \textit{hole} is defined as a connected region of local intensity minima (dark pixels) that are completely surrounded by pixels of higher intensity. The parameter \textit{conn} specifies the connectivity (e.g. $4$ or $8$, see Figure \ref{conn}) that determines how neighboring pixels are considered to be part of the same region.\\
The algorithm used by \texttt{imfill} is based on the concept of geodesic dilation morphological reconstruction \cite{imfill}.\\
\begin{figure}[htbp]
    \centering
    \begin{subfigure}{0.1\textwidth}
        \includegraphics[width=\linewidth]{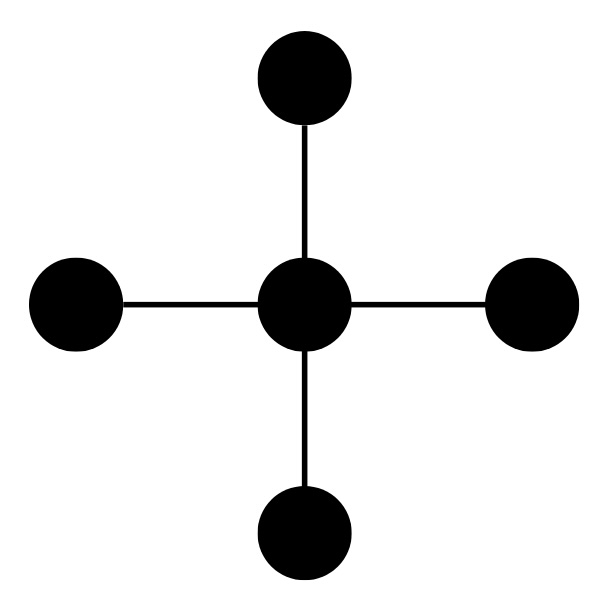}
        \caption{}
    \end{subfigure}
    \hspace{0.05\textwidth}
    \begin{subfigure}{0.1\textwidth}
        \includegraphics[width=\linewidth]{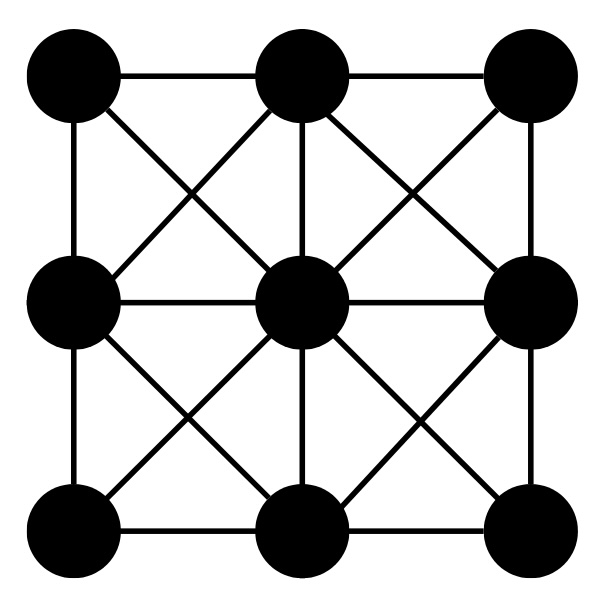}
        \caption{}
    \end{subfigure}
\caption{The elementary ball in $4$-connectivity (a), and $8$-connectivity (b).}
\label{conn}
\end{figure}

We also evaluate our method against two additional inpainting techniques available in \textsc{Matlab\textsuperscript{\textcopyright}}: {\tt{inpaintExemplar(A,mask)}} and {\tt{inpaintCoherent(A,mask)}}, where \textit{mask} is a logical image indicating the target regions to be inpainted. The coherent method \cite{inpaintCoherent} is a diffusion-based approach that propagates image structures or level lines via diffusion based on partial differential equations and variational methods. While this method is efficient and fast, it may introduce blurring in large gap regions. The exemplar-based method \cite{inpaintExemplar1, inpaintExemplar2} operates in a patch-wise manner. It computes patch priorities at the boundary of the missing region and fills in the target region by copying best-matching patches from the surrounding (known) image area. This approach is particularly effective for textured regions and preserves structural information well, although it can be more computationally demanding.\\

Finally, we include a comparison with the Mixed Median inpainting method \cite{mixedMedian}. This approach is a non-iterative, filter-based technique designed to restore missing or corrupted regions by exploiting robust statistical information from the local neighborhood. In particular, the method replaces missing pixels by applying a mixed median operator that combines directional and spatial median filtering strategies, allowing for more stable reconstruction in the presence of noise and structural discontinuities. Unlike diffusion-based approaches, the mixed median method does not rely on solving partial differential equations, and unlike patch-based techniques, it does not require explicit patch search or matching procedures. This makes it computationally efficient while still preserving edge information reasonably well. However, its performance may degrade in large missing regions or highly textured areas, where the lack of long-range structural modeling can lead to over-smoothing or loss of fine details.

\begin{figure}[htbp]
    \centering
    \begin{subfigure}{0.22\textwidth}
        \includegraphics[width=\linewidth]{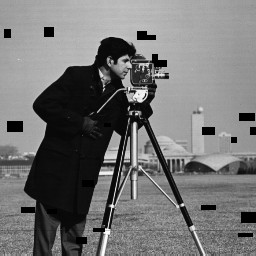}
        \caption{}
    \end{subfigure}
    \begin{subfigure}{0.22\textwidth}
        \includegraphics[width=\linewidth]{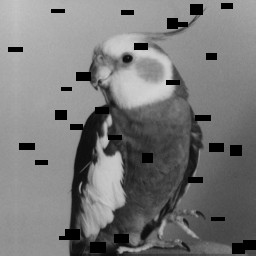}
        \caption{}
    \end{subfigure}
    \begin{subfigure}{0.22\textwidth}
        \includegraphics[width=\linewidth]{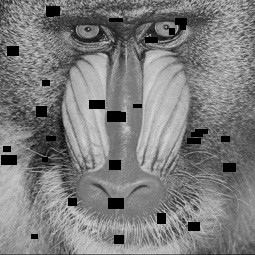}
        \caption{}
    \end{subfigure}
\caption{Images with simulated loss: (a) \textit{Cameraman} with $2.43$\% missing pixels; (b) \textit{Bird} with $4.47$\% missing pixels; (c) \textit{Baboon} with $3.79$\% missing pixels.}
\label{images_with_gaps}
\end{figure}

To test our method, we select three grayscale images: \textit{Cameraman}, \textit{Bird}, and \textit{Baboon}. These images are part of the Grayscale Set 1 and 2, provided by the Waterloo Fractal Coding and Analysis Group, and are available in the repository at \cite{url3}. Artificial black pixels are introduced (see Figure \ref{images_with_gaps}). Since the true pixel values are known, we directly evaluate the prediction accuracy using the well-known Structural Similarity Index (SSIM) and Peak Signal-to-Noise Ratio (PSNR); their formulations are provided in Table \ref{Indici2}.
\begin{table}[tbph]
\centering
\renewcommand{\arraystretch}{1.5} % Aumenta l'interlinea
\begin{tabular}{ll}
\hline
\rowcolor[HTML]{EFEFEF} 
\textbf{\small Gap-filling indexes}                           & \multicolumn{1}{c}{\cellcolor[HTML]{EFEFEF}} \\
\small mean square error & {\footnotesize $\displaystyle\text{MSE}:=\frac{1}{nm}\sum_{i=1}^{n}\sum_{j=1}^{m}|A(i,j)-B(i,j)|^2$} \\
\rowcolor[HTML]{EFEFEF} 
\small peak signal-to-noise ratio & {\footnotesize $\displaystyle\text{PSNR}:=20\log_{10}\left(\frac{\max{A}}{\sqrt{\text{MSE}}}\right)$} \\
\small structural similarity index & {\footnotesize $\displaystyle\text{SSIM}:=\frac{2\mu_{A} \mu_{B}\, + C_1}{\mu_{A}^2 + \mu_{B}^2\, + C_1}\cdot\frac{2\sigma_{AB}\, + C_2}{\sigma_{A}^2 + \sigma_{B}^2 + C_2}$} \\ 
\hline
\end{tabular}
\caption{Indexes for evaluating gap-filling performance. In the formulas, $A$ denotes the original image, while $B$ represents the reconstructed image obtained after applying the gap-filling techniques.}
\label{Indici2}
\end{table}

Figure \ref{images_with_SKprediction} shows the reconstructed images obtained using LP-SK algorithm with $w=40$ and $s=9$. However, similar results were already obtained with $w>20$ and $s>5$, which means the method performs well even with less demanding parameters. Additionally, the values of $w$ and $s$ do not affect the computational cost (i.e., the CPU time), as the algorithm runs in very small fractions of a second regardless of these settings. In Figure \ref{plots}, we display the absolute error for each reconstructed pixel. The error values are generally close to zero, with some peaks occurring in structurally complex or non-homogeneous regions, where reconstruction is more challenging.\\
\begin{figure}[htbp]
    \centering
    \begin{subfigure}{0.22\textwidth}
        \includegraphics[width=\linewidth]{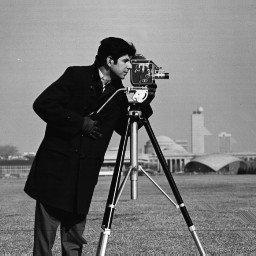}
        \caption{}
    \end{subfigure}
    \begin{subfigure}{0.22\textwidth}
        \includegraphics[width=\linewidth]{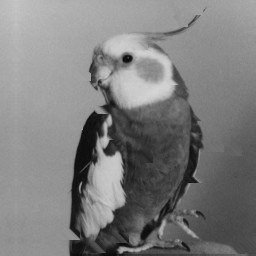}
        \caption{}
    \end{subfigure}
    \begin{subfigure}{0.22\textwidth}
        \includegraphics[width=\linewidth]{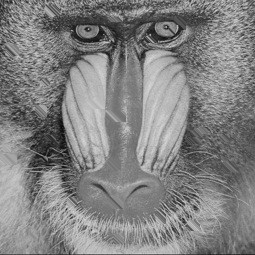}
        \caption{}
    \end{subfigure}
\caption{Reconstructed images using LP-SK Algorithm.}
\label{images_with_SKprediction}
\end{figure}
\begin{figure}[htbp]
    \centering
    \begin{subfigure}{0.32\textwidth}
        \includegraphics[width=\linewidth]{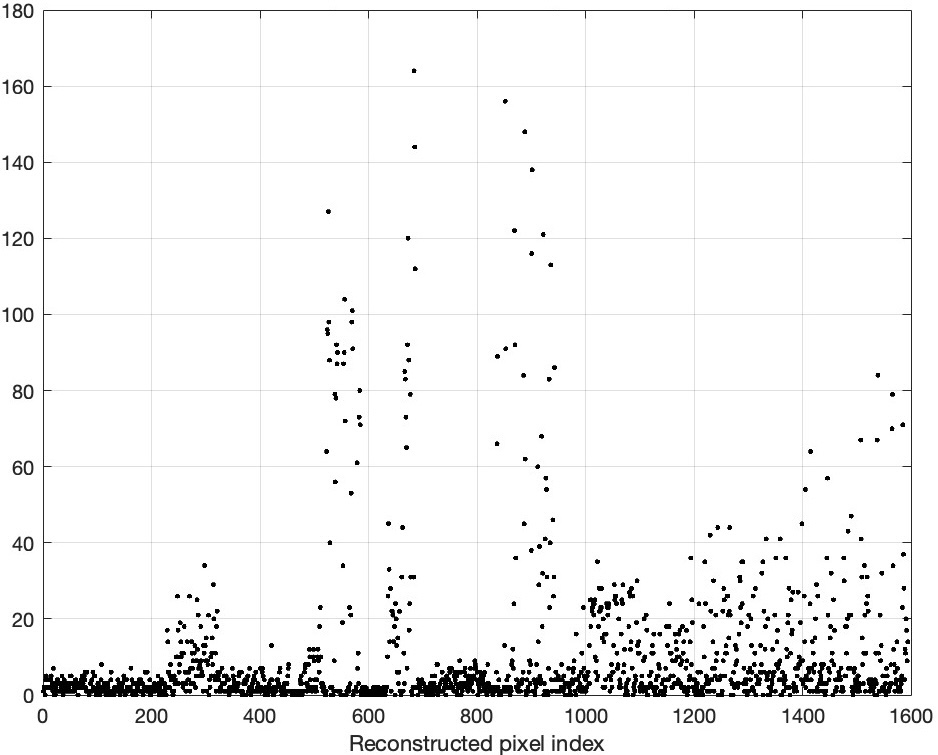}
        \caption{}
    \end{subfigure}
    \begin{subfigure}{0.32\textwidth}
        \includegraphics[width=\linewidth]{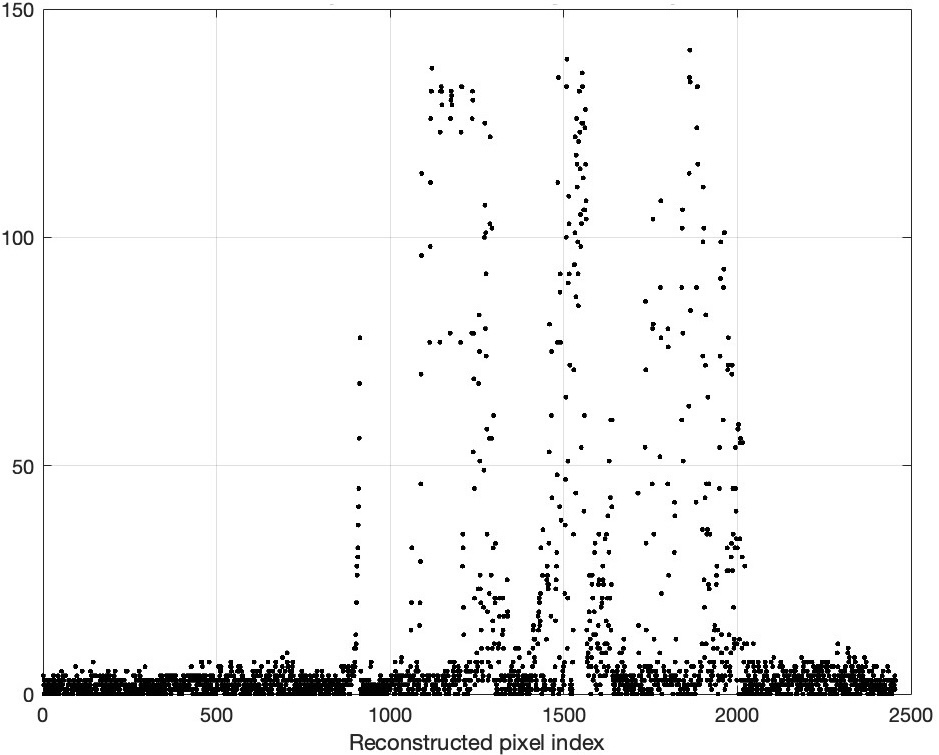}
        \caption{}
    \end{subfigure}
    \begin{subfigure}{0.32\textwidth}
        \includegraphics[width=\linewidth]{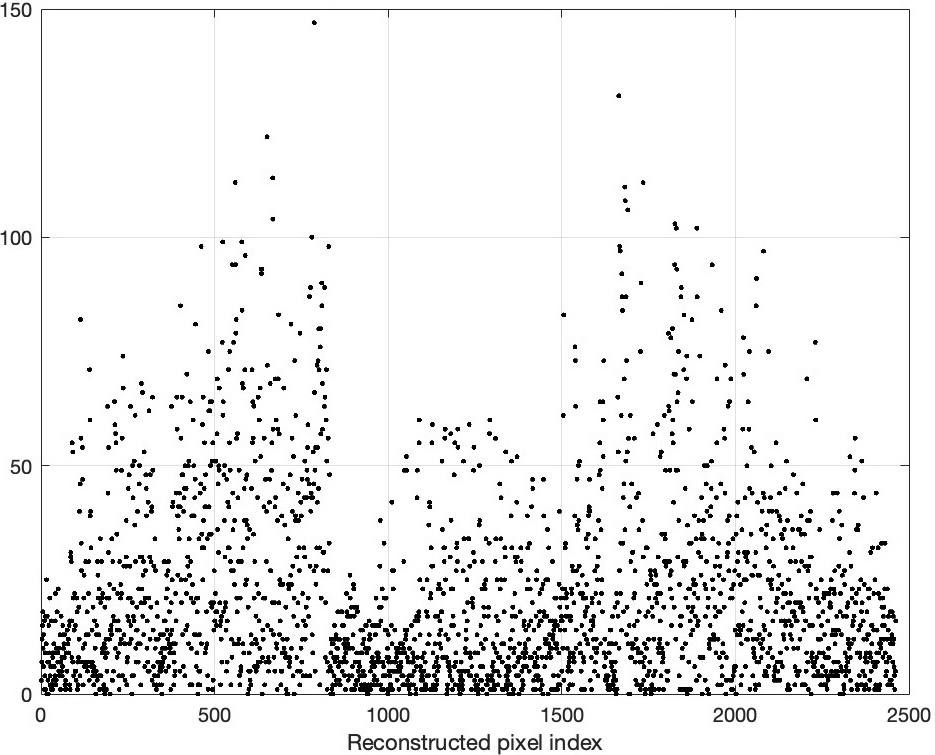}
        \caption{}
    \end{subfigure}
\caption{Absolute error for each pixel reconstructed by LP-SK: (a) \textit{Cameraman}; (b) \textit{Bird}; (c) \textit{Baboon}.}
\label{plots}
\end{figure}

Figure \ref{image_text} presents an example in which we added synthetic text to the \textit{Cameraman} image, covering $11.44$\% of its area, to simulate an extensive missing data scenario. \\
\begin{figure}[htbp]
    \centering
    \begin{subfigure}{0.22\textwidth}
        \includegraphics[width=\linewidth]{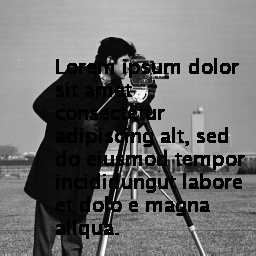}
        \caption{}
    \end{subfigure}
    \begin{subfigure}{0.22\textwidth}
        \includegraphics[width=\linewidth]{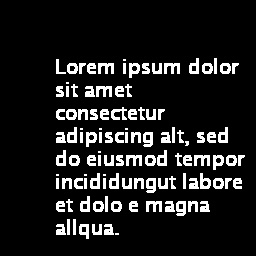}
        \caption{}
    \end{subfigure}
    \begin{subfigure}{0.22\textwidth}
        \includegraphics[width=\linewidth]{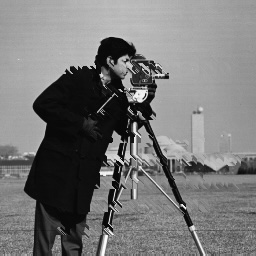}
        \caption{}
    \end{subfigure}
\caption{Removal of superimposed text: (a) \textit{Cameraman} ($11.44$\% covered by text); (b) corresponding mask; (c) reconstructed image using LP-SK.}
\label{image_text}
\end{figure}

Quantitative results are shown in Table \ref{tab_gap_filling}, where we report the SSIM and PSNR values obtained with LP-SK algorithm, \texttt{imfill}, \texttt{inpaintExemplar}, \texttt{inpaintCoherent}, and \texttt{mixedMedian} . Among these methods, \texttt{inpaintCoherent} consistently achieves the highest performance. On the other hand, the \texttt{imfill} method shows the weakest results, as it not only fails to reconstruct the missing pixels correctly but also alters other parts of the image, in particular, it tends to modify all dark pixels surrounded by brighter ones. This is visually illustrated in Figure \ref{other_images_reconstructed}(a). Missing pixels that are located in already dark areas of the image are often not reconstructed, while valid dark areas, such as the eye in the \textit{Cameraman} image, are incorrectly filled, as if they were holes. Figures \ref{other_images_reconstructed}(b) and \ref{other_images_reconstructed}(c) display examples of the application of \texttt{inpaintExemplar} on the \textit{Bird} image and \texttt{inpaintCoherent} on the \textit{Baboon} image, respectively. Figure \ref{other_images_reconstructed}(d) shows the \textit{Cameraman} image with superimposed text reconstructed using \texttt{mixedMedian}, illustrating its ability to recover local structures while preserving overall visual coherence, although fine details in the corrupted regions may be partially smoothed due to its local median-based formulation.\\

\begin{table}[htbp]
\centering
\resizebox{\columnwidth}{!}{%
\begin{tabular}{lccccc c}
\hline
                                    &      & LP-SK algorithm & imfill  & inpaintExemplar & inpaintCoherent & mixedMedian \\ \hline
\multirow{2}{*}{\textit{Cameraman}} & PSNR & 37.7977              & 25.2094 & 38.0534          & 39.2661 & 38.4096 \\
                                    & SSIM & 0.9902               & 0.8963  & 0.9890           & 0.9936 & 0.9924 \\ \hline
\multirow{2}{*}{\textit{\begin{tabular}[c]{@{}l@{}}Cameraman\\ covered by text\end{tabular}}} & PSNR & 23.6036 & 20.8821 & 26.3141 & 28.4632 & 26.6339 \\
                                 & SSIM & 0.9241               & 0.8274  & 0.9403           & 0.9648  & 0.9521 \\ \hline
\multirow{2}{*}{\textit{Bird}}      & PSNR & 32.8068              & 26.7577 & 34.3324          & 40.7840 & 38.9450 \\
                                    & SSIM & 0.9802               & 0.9167  & 0.9814           & 0.9931  & 0.9910 \\ \hline
\multirow{2}{*}{\textit{Baboon}}    & PSNR & 32.7452              & 23.5143 & 33.1104          & 35.0852 & 35.4458 \\
                                    & SSIM & 0.9705               & 0.7998  & 0.9707           & 0.9783 & 0.9788 \\ \hline
\end{tabular}%
}
\caption{Performance comparison of gap filling methods. These results refer to the complete images.}
\label{tab_gap_filling}
\end{table}

Although \texttt{inpaintCoherent} currently outperforms LP-SK in terms of reconstruction quality, the proposed approach shows potential as a fast and efficient solution for gap filling. One key advantage of the LP-SK approach is that it works using only past pixel values, without needing information from the rest of the image. This makes it especially suitable for applications where data become available over time, such as in satellite data streams or real-time remote sensing, where future or complete image data may not yet be available. In contrast, inpainting algorithms typically use information from pixels that is not necessarily spatially preceding the missing pixels. Given its simplicity and efficiency, LP-SK offers a good starting point, and we plan to handle more complex scenarios or larger gaps. 
\begin{figure}[htbp]
    \centering
    \begin{subfigure}{0.22\textwidth}
        \includegraphics[width=\linewidth]{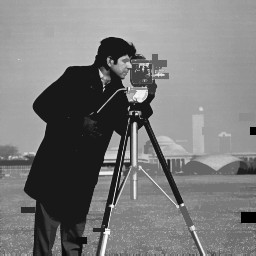}
        \caption{}
    \end{subfigure}
    \begin{subfigure}{0.22\textwidth}
        \includegraphics[width=\linewidth]{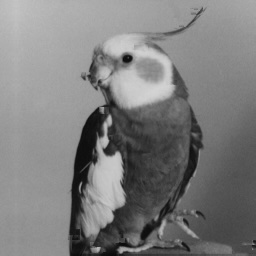}
        \caption{}
    \end{subfigure}
    \begin{subfigure}{0.22\textwidth}
        \includegraphics[width=\linewidth]{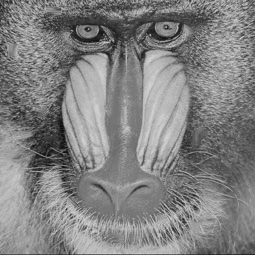}
        \caption{}
    \end{subfigure}
    \begin{subfigure}{0.22\textwidth}
        \includegraphics[width=\linewidth]{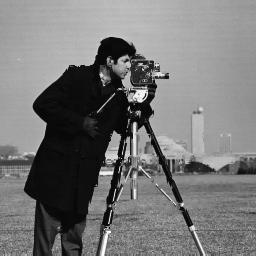}
        \caption{}
    \end{subfigure}
\caption{Reconstructed images using other gap filling methods: (a) \textit{Cameraman} with {\tt{imfill}}; (b) \textit{Bird} with {\tt{inpaintExemplar}}; (c) \textit{Baboon} with {\tt{inpaintCoherent}}; (d) \textit{Cameraman} with superimposed text reconstructed using {\tt{mixedMedian}}.}
\label{other_images_reconstructed}
\end{figure}
In future work, we also plan to extend the comparison of the proposed LP-SK method to more recent inpainting approaches, in particular deep learning-based methods, which have shown strong performance in the literature. However, these techniques rely on training data and model design, which is beyond the scope of this work focused on classical prediction-based methods. For a broader overview of recent developments in image inpainting, we refer the reader to \cite{deepinpainting_survey1, deepinpainting_survey2}.
\section{Speckle noise and its reduction}\label{s4}
\textit{Speckle noise} is a granular disturbance that affects SAR images and all coherent imaging techniques, such as ultrasound and laser imaging. Unlike Gaussian or salt-and-pepper noise, which are additive, speckle noise is multiplicative in nature.
Let $A_n(i, j)$ denote a distorted pixel in an image, and let $M(i, j)$ denote the corresponding noise-free image pixel. According to the multiplicative noise model,
\begin{equation*}\label{errore_moltiplicativo}
A_n(i, j) = M(i, j)\cdot N(i, j),
\end{equation*}
where $N(i, j)$ is independent of $M(i, j)$ and represents the speckle fading term, whose mean value is $1$. Finally, $(i,j)$ represents the spatial position of the pixels. This means that it affects the variance of the pixel values rather than their mean intensity. The greater the variance of speckle noise, the more difficult it is to recover a noiseless image. This results in a textured appearance that can obscure fine details, reduce contrast, and distort gray values. To better explain the influence of speckle noise, Figure \ref{speckle1} presents an image with a gray value of 40, affected by speckle noise with variances of 0.1 and 0.9. Further, Figure \ref{histo} shows the histograms of the original image and the noisy images with variances of 0.05 and 0.2. As noise increases, the histogram gradually shifts towards a more uniform distribution, indicating a gradual equalization effect.
\begin{figure}[tbph] 
\begin{center}
\includegraphics[width=13cm]{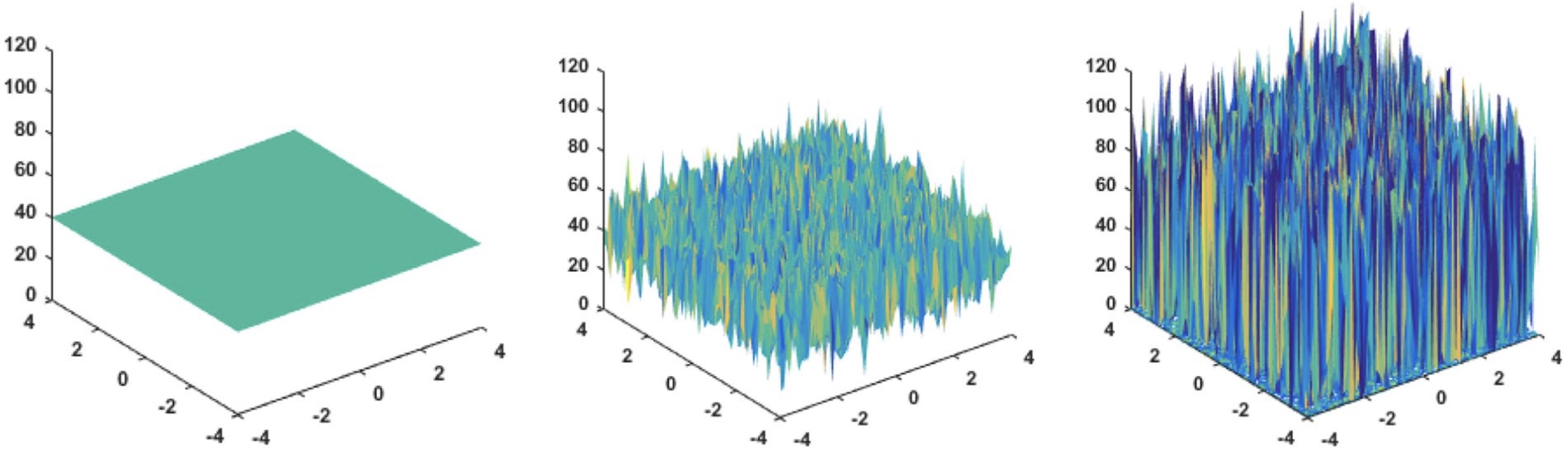}
\end{center}
\caption{Effects of speckle noise with variances of 0.1 and 0.9 on image gray values. Reproduced from \cite{speckle2} with permission from Springer Nature.}
\label{speckle1}
\end{figure}
\begin{figure}[htbp]
    \centering
    \begin{subfigure}{0.3\textwidth}
        \includegraphics[width=\linewidth]{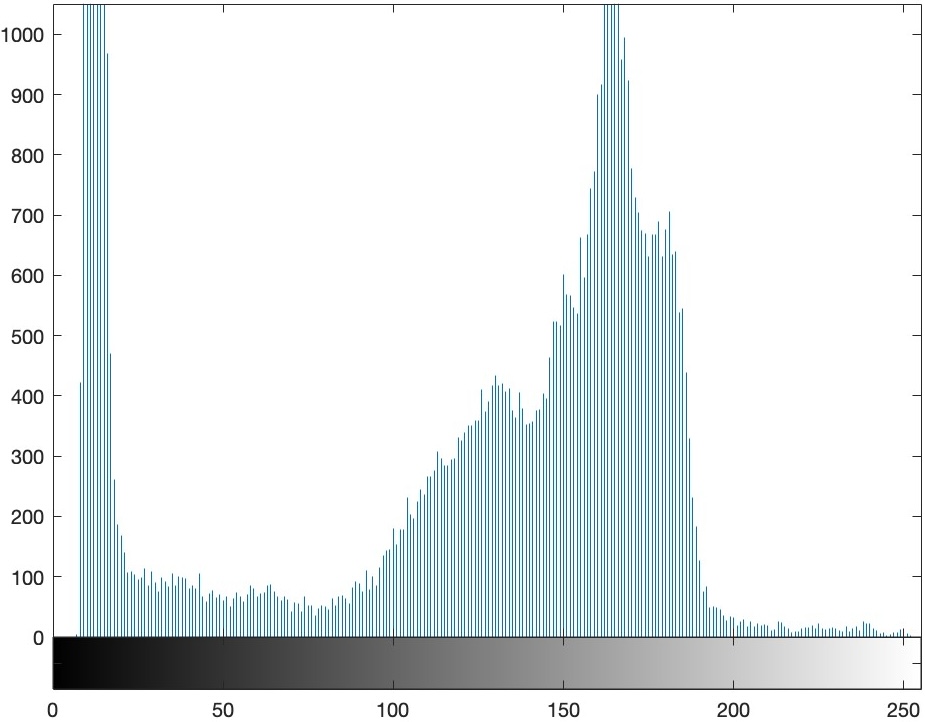}
        \caption{}
    \end{subfigure}
    \begin{subfigure}{0.3\textwidth}
        \includegraphics[width=\linewidth]{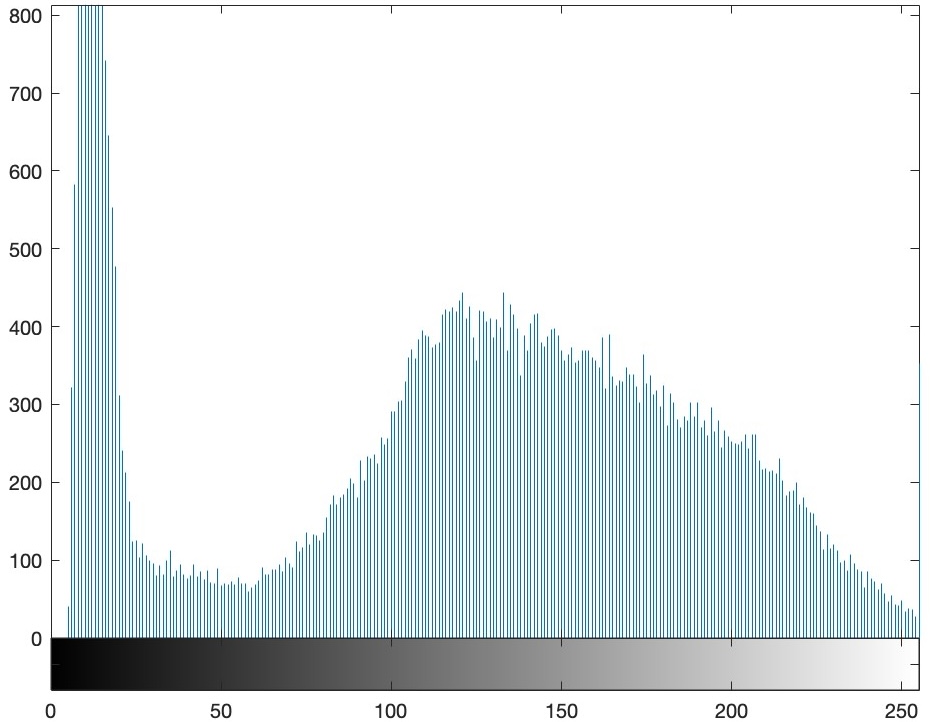}
        \caption{}
    \end{subfigure}
    \begin{subfigure}{0.3\textwidth}
        \includegraphics[width=\linewidth]{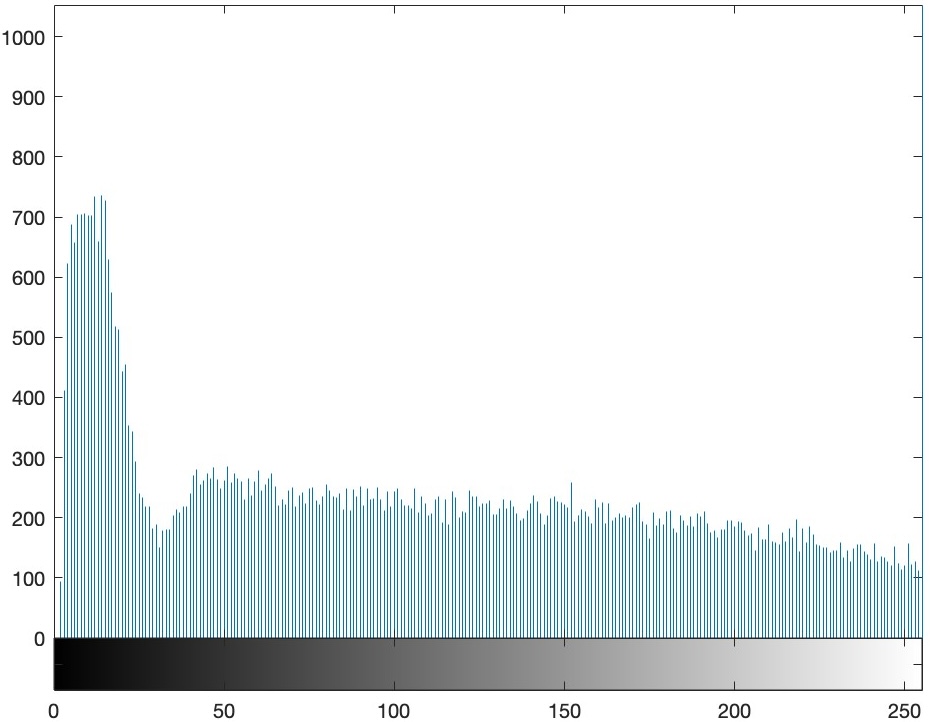}
        \caption{}
    \end{subfigure}
\caption{Histograms of \textit{Cameraman} image \cite{url3}: (a) original; (b) speckled with $\sigma^2=0.05$; (c) speckled with $\sigma^2=0.2$.}
\label{histo}
\end{figure}
\\

During the last three decades, several methods have been proposed for the reduction of speckle, or \textit{despeckling}, in SAR images \cite{despeckle1}. The existence of speckle noise in SAR images is an inherent and specific characteristic which is deterministic and random in
nature. Classical algorithms include common mean filter, median filter, Lee and Frost filters, and non-local mean filter. 
\begin{itemize}
    \item \textbf{Mean filter} (also known as the \textit{boxcar filter}) is an averaging filter that replaces the center pixel in a $3\times 3$ or a larger moving kernel $n\times m$ with the mean value of the surrounding pixels. It can be implemented using the following formula: 
    \begin{equation*}
    A_d(i,j):=\frac{1}{nm}\sum_{(k,l)\in S}A_n(k,l),
    \end{equation*}
    where $S$ is the $n\times m$ neighborhood of the pixel, $I_d (i,j)$ is the processed image and $A_n(k,l)$ is the input image. This filter is simple and fast; however, it also causes a reduction in detail and resolution.
    \item \textbf{Median filter} is a non-linear filter that replaces the center pixel in a kernel with the median value of its neighbors. It is particularly effective, as it removes speckle noise without excessively blurring edges, making it well-suited for images with sharp transitions and textures. Compared to the mean filter, it preserves edges better but may distort very small details if speckle noise is present. Since the median is less sensitive than the mean to extreme values, it removes outliers more effectively than the mean filter.
    \item \textbf{Frost filter} \cite{Frost} replaces the pixel of interest with a weighted sum of the values within an $n\times n$ kernel. In particular, the filtered value is a linear combination of pixel values within the local window with a Gaussian weighting function that depends on the local coefficient of variation $\overline{\sigma}$ of the noisy image $I_n$, and on the ratio of local standard deviation $\sigma_{I_n}$ to local mean $\mu_{I_n}$:
    \begin{equation*}
    \text{Digital number (DN)}=\sum_{n\times n} k \alpha e^{-\alpha |t|},
    \end{equation*}
    where $k$ is a normalization constant, $\alpha:=(4/n\overline{\sigma}^2)\cdot(\sigma_{I_n}^2/\mu^2_{I_n})$, and $|t|=|x-x_0|+|y-y_0|$.\\
    In uniform regions, the Frost filter behaves like a mean filter, smoothing out speckle. In high-contrast regions, it acts as a high-pass filter, with rapid decay of elements away from the filter center. Thus, uniform areas are smoothed, removing speckle, while high-contrast edges and objects retain their signal values without smoothing.
    \item \textbf{Lee filter} \cite{Lee} is considered the first model-based despeckling filter. This filter is adaptive to the local statistics in an image and it is based on the minimum mean square error (MMSE). The despeckled image $A_d$ is produced on the following equation:
    \begin{equation*}
        A_d= \text{mean}+K \left(A_n - \text{mean}\right),
    \end{equation*}
    where ``mean'' is the avarage of pixels in a moving kernel, $K:= \text{var}(x)/(\text{mean}^2\sigma_{I_n}^2+\text{var}(x))$, $\text{var}(x):=(\text{variance within kernel} + [\text{mean within kernel}]^2)/(\sigma_{I_n}^2+1)-[\text{mean within kernel}]^2$. Lee filter is an isotropic adaptive filter which can not remove noise in the edge region effectively. Despite this, it is good at preserving prominent edges, linear features, point target, and texture information.
    \item \textbf{Non-local mean (NLM) filter} \cite{NLM} works by estimating a noise-free image as a weighted average of noisy pixels. It performs well with additive white Gaussian noise, however, in the case of SAR images, the weights have to be generalized to the case of multiplicative and non-Gaussian speckle.
\end{itemize}
\subsection{Evaluation metrics for despeckling efficiency}\label{ss4_1}
Subjective evaluations concern how the textural details and visual effect of an image are improved after denoising. Objective evaluations, on the other hand, rely on quantitative metrics; see Table \ref{Indici}. Some of these metrics, namely MSE, PSNR, and SSIM, were previously introduced in the gap-filling framework (see Table \ref{Indici2}) and are reported again here for completeness.\\
The notations in Table \ref{Indici} are as follows: $A_{n}$ denotes the noisy image, $A_d$ represents the despeckled image, and $A$ is the ground-truth noiseless image. The image dimensions are denoted by $n$ and $m$. Metrics that require $A$ are referred to as \textit{with-reference metrics}, whereas those that do not depend on $A$ are classified as \textit{without-reference metrics}.
\begin{table}[tbph]
\centering
\renewcommand{\arraystretch}{1.5} % Aumenta l'interlinea
\begin{tabular}{ll}
\hline
\rowcolor[HTML]{EFEFEF} 
\textbf{\small With-reference indexes}                           & \multicolumn{1}{c}{\cellcolor[HTML]{EFEFEF}} \\
\small mean square error                                   & {\footnotesize $\displaystyle\text{MSE}:=\frac{1}{nm}\sum_{i=1}^{n}\sum_{j=1}^{m}|A(i,j)-A_d(i,j)|^2$}                                      \\
\rowcolor[HTML]{EFEFEF} 
\small peak signal-to-noise ratio                         & {\footnotesize $\displaystyle\text{PSNR}:=20\log_{10}\left(\frac{\max{A}}{\sqrt{\text{MSE}}}\right)$}                                      \\
\small structural similarity index measurement      & {\footnotesize $\displaystyle\text{SSIM}:=\frac{2\mu_{A} \mu_{A_d}\, + C_1}{\mu_{A}^2 + \mu_{A_d}^2\, + C_1}\cdot\frac{2\sigma_{AA_d}\, + C_2}{\sigma_{A}^2 + \sigma_{A_d}^2 + C_2}$}                                      \\ \hline
\rowcolor[HTML]{EFEFEF} 
{\color[HTML]{000000} \textbf{\small Without-reference indexes}} &      
\\
\small speckle intensity                           & {\footnotesize $\displaystyle\text{SI}:= \frac{\sqrt{\sigma_{A_d}}}{\mu_{A_d}}$ }                                     \\
\rowcolor[HTML]{EFEFEF} 
\small speckle suppression index                           & {\footnotesize $\displaystyle\text{SSI}:= \frac{\sqrt{\sigma_{A_d}}}{\mu_{A_d}}\cdot \frac{\mu_{A_n}}{\sqrt{\sigma_{A_n}}}$ }                                      \\
\small speckle mean preservation index                                    & {\footnotesize $\displaystyle\text{SMPI}:=(1+|\mu_{A_n}-\mu_{A_d}|)\cdot \sqrt{\frac{\sigma_{A_d}}{\sigma_{A_n}}}$ }                              
\\ \rowcolor[HTML]{EFEFEF} 
\small equivalent number of looks                                    & {\footnotesize $\displaystyle\text{ENL}:=\left(\frac{\mu_{A_d}}{\sigma_{A_d}}\right)^2$}                                      \\ \hline
\end{tabular}
\caption{Indexes for evaluating despeckling performance.}
\label{Indici}
\end{table}

We stress that without-reference indexes do not rely on the knowledge of the ground-truth. Instead, they are based only on specific statistical hypotheses about the signal model. Since the signal model is strongly dependent on the degree of scene heterogeneity, it is necessary to select appropriate areas, e.g. homogeneous regions, for the computation of such specific indexes.\\

In particular, the \textit{speckle suppression index} (SSI) \cite{SX1996} is obtained by normalizing the \textit{speckle intensity} (SI) of the despeckled image \cite{SIindex}, and the SI of the original noisy image in a specific homogeneous area. The lower the SSI is ($SSI\ll 1$), the better despeckling performance the filter has. However, SSI may not evaluate speckle removal performance if the filter overestimates the mean of the speckled image. The \textit{speckle mean preservation index} (SMPI) was introduced in \cite{WangGe}, and lower SMPI values show better filter efficiency. The \textit{equivalent number of looks} (ENL) \cite{Bruniquel} is also known as measure of the signal-to-noise
ratio. The higher the ENL value for a filter, the higher the efficiency of improving the signal-to-speckle ratio over homogeneous areas so that mean image quality is better.
\subsection{Proposed Down-Up processing algorithm}\label{ss4_2}
This section introduces the Down-Up processing algorithm aimed at reducing speckle noise in digital images. The process starts with a noise-free grayscale image of size $n\times m$, which serves as a reference. To simulate real conditions, we add speckle noise with a variance of $0.05$ using \textsc{Matlab\textsuperscript{\textcopyright}} inbuilt command function {\tt{imnoise(A,"speckle",0.05)}}, producing the noisy image. Such image is then processed using the five filters described previously. To evaluate the filtering performance, we compute both with-reference and without-reference quality metrics, with the latter measured over two homogeneous regions of interest (ROIs). These ROIs have been selected in non-edge areas of the images in order to avoid boundary effects and ensure a consistent statistical evaluation.

In addition to the standard filtering approach, we propose an alternative method that introduces two additional steps: downscaling and upscaling. For this reason, we refer to it as the \textit{Down-Up processing algorithm} (see Figure \ref{schema}). First, the noisy image is reduced to dimensions $\frac{n}{2}\times \frac{m}{2}$ using three different downscaling techniques: \textit{SK algorithm}, \textit{bilinear} and \textit{bicubic interpolation} (see, e.g., \cite{keys1981cubic, smith1981bilinear, getreuer2011linear}). These smaller images, still affected by noise, are then processed using the same five filters. After filtering, we restore the images to their original size using the same three upscaling methods. Once back to their original resolution, the quality metrics are computed and compared with those obtained from direct filtering. The aim is to show that this Down-Up approach improves noise reduction.
\begin{figure}[tbph] 
\begin{center}
\includegraphics[width=\textwidth]{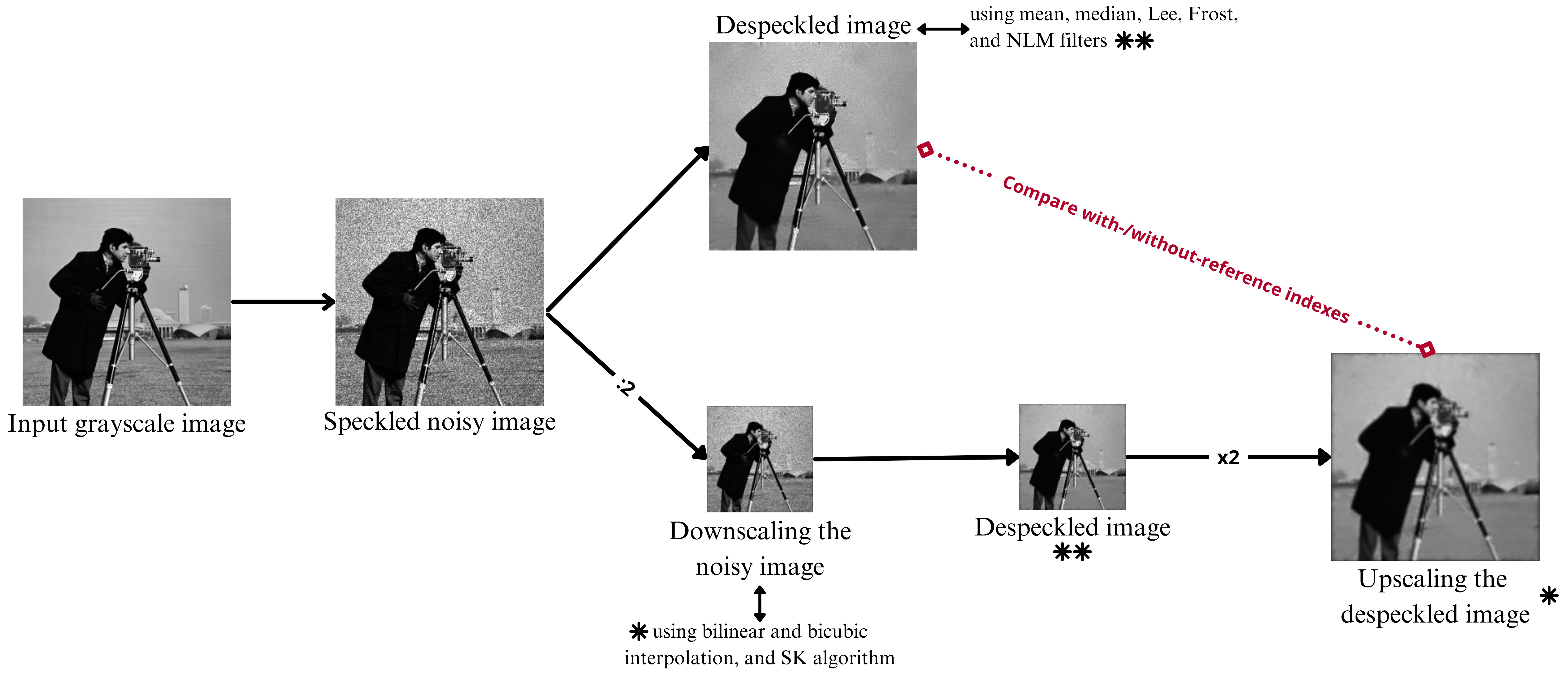}
\end{center}
\caption{Scheme of the proposed Down-Up scaling method to enhance speckle reduction efficiency.}
\label{schema}
\end{figure}

\subsection{Numerical testing}\label{ss4_3}
In the present section, both simulated and real SAR images are used to analyze the experimental results in terms of quantitative evaluations, assessing the performance of speckle noise removal.\\

First, we select from \cite{url3} two 8-bit grayscale images \textit{Cameraman} with resolution $256\times 256$, and \textit{Boat} with resolution $512\times 512$ (see, Figure \ref{cameraman_roi}).\\
\begin{figure}[tbph] 
\centering
\begin{subfigure}{0.22\textwidth}
        \includegraphics[width=\linewidth]{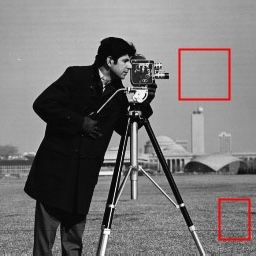}
        \caption{}
\end{subfigure}
\begin{subfigure}{0.22\textwidth}
        \includegraphics[width=\linewidth]{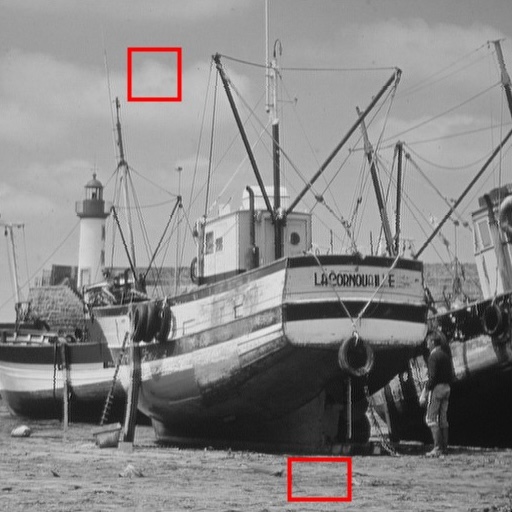}
        \caption{}  
\end{subfigure}
\caption{(a) \textit{Cameraman}: ROI1$=[220,200,30,40]$, and ROI2$=[180,50,50,50]$; (b) \textit{Boat}: ROI1$=[130,50,50,50]$, and ROI2$=[290,460,60,40]$.}
\label{cameraman_roi}
\end{figure}

The results of our experiments are detailed in Tables \ref{tab_cameraman} and \ref{tab_boat}, which summarize the quantitative with- and without-reference indexes for different Down-Up configurations applied to \textit{Cameraman} and \textit{Boat}, respectively. The best-performing metrics are highlighted in bold for easy identification. The configurations involve different combinations of SK algorithm, bilinear and bicubic interpolation applied before and after filtering. Our experiments show a general improvement in without-reference quality metrics. However, the highest PSNR values are still obtained when filters are applied directly to the noisy image. Further, these results align with previous studies \cite{CSV, CAR}, which suggest that SK is better for upscaling, while bicubic interpolation seems to be more effective for downscaling. Based on these findings, the best approach seems to be bicubic downscaling, followed by filtering and SK upscaling. Although this approach does not maximize SSIM and PSNR, it significantly improves without-reference metrics, which are the most useful in real-world applications where noise-free reference images are usually unavailable. To reinforce these findings, we also tested other additional grayscale images. However, to keep the paper concise, we have not included the numerical results. Nonetheless, the observed trends align with those found for \textit{Cameraman} and \textit{Boat}.\\ 
\begin{table}[tbph]
\centering
\resizebox{\columnwidth}{!}{%
\begin{tabular}{lcccccccccc}
\hline
\multicolumn{1}{c}{} &
   &
   &
  \multicolumn{2}{c}{SI} &
  \multicolumn{2}{c}{SSI} &
  \multicolumn{2}{c}{SMPI} &
  \multicolumn{2}{c}{ENL} \\ \cline{4-11} 
\multicolumn{1}{c}{\multirow{-2}{*}{}} &
  \multirow{-2}{*}{PSNR} &
  \multirow{-2}{*}{SSIM} &
  ROI1 &
  ROI2 &
  ROI1 &
  ROI2 &
  ROI1 &
  ROI2 &
  ROI1 &
  ROI2 \\ \hline
mean           & \textbf{23.8680} & 0.5726          & 0.4660 & 0.3578 & 0.6299 & 0.6060 & 0.6294          & 0.6059          & 99.6423   & 147.8924 \\
\rowcolor[HTML]{EFEFEF} 
SK+mean+SK     & 21.0783          & 0.6129          & 0.3254 & 0.2638 & 0.4399 & 0.4469 & \textbf{0.4409} & 0.4463          & 416.1002  & 503.1320 \\
bil+mean+bil   & 21.7487          & 0.6382          & 0.3328 & 0.2612 & 0.4498 & 0.4424 & 0.4648          & 0.4589          & 365.0157  & 497.5693 \\
\rowcolor[HTML]{EFEFEF} 
bil+mean+SK    & 21.4919          & 0.6230          & 0.3188 & 0.2550 & 0.4309 & 0.4320 & 0.4449          & 0.4478          & 433.8961  & 547.6087 \\
bic+mean+bic   & 21.6372          & \textbf{0.6424} & 0.3283 & 0.2586 & 0.4437 & 0.4381 & 0.4608          & 0.4559          & 382.9099  & 515.4174 \\
\rowcolor[HTML]{EFEFEF} 
bic+mean+SK &
  21.3618 &
  0.6227 &
  \textbf{0.3145} &
  \textbf{0.2528} &
  \textbf{0.4250} &
  \textbf{0.4282} &
  0.4411 &
  \textbf{0.4454} &
  \textbf{455.2192} &
  \textbf{564.8304} \\ \hline
median         & \textbf{22.2473} & 0.4838          & 0.5421 & 0.4352 & 0.7328 & 0.7371 & 0.7292          & 0.7400          & 55.2257   & 67.1847  \\
\rowcolor[HTML]{EFEFEF} 
SK+median+SK   & 21.1981          & 0.6211          & 0.3417 & 0.2626 & 0.4619 & 0.4447 & 0.4604          & \textbf{0.4445} & 347.5583  & 510.5296 \\
bil+median+bil & 22.0316          & 0.6363          & 0.3500 & 0.2706 & 0.4731 & 0.4583 & 0.4862          & 0.4798          & 300.4479  & 427.1522 \\
\rowcolor[HTML]{EFEFEF} 
bil+median+SK  & 21.7365          & 0.6325          & 0.3317 & 0.2623 & 0.4484 & 0.4443 & 0.4601          & 0.4647          & 373.2805  & 484.2200 \\
bic+median+bic & 21.6372          & \textbf{0.6424} & 0.3283 & 0.2652 & 0.4437 & 0.4491 & 0.4608          & 0.4709          & 382.9099  & 462.1961 \\
\rowcolor[HTML]{EFEFEF} 
bic+median+SK &
  21.6817 &
  0.6353 &
  \textbf{0.3255} &
  \textbf{0.2579} &
  \textbf{0.4400} &
  \textbf{0.4369} &
  \textbf{0.4537} &
  0.4578 &
  \textbf{399.7759} &
  \textbf{516.7099} \\ \hline
Frost          & \textbf{23.9332} & 0.5739          & 0.4665 & 0.3579 & 0.6307 & 0.6061 & 0.6302          & 0.6060          & 99.1377   & 147.8118 \\
\rowcolor[HTML]{EFEFEF} 
SK+Frost+SK    & 20.6573          & 0.6120          & 0.3253 & 0.2639 & 0.4396 & 0.4470 & \textbf{0.4406} & 0.4464          & 417.0025  & 502.6994 \\
bil+Frost+bil  & 21.6573          & 0.6208          & 0.3328 & 0.2612 & 0.4498 & 0.4424 & 0.4647          & 0.4589          & 365.1845  & 497.6040 \\
\rowcolor[HTML]{EFEFEF} 
bil+Frost+SK   & 21.1727          & \textbf{0.6211} & 0.3189 & 0.2551 & 0.4310 & 0.4321 & 0.4449          & 0.4479          & 433.5581  & 547.1650 \\
bic+Frost+bic  & 21.4809          & \textbf{0.6211} & 0.3285 & 0.2586 & 0.4440 & 0.4380 & 0.4611          & 0.4558          & 381.9714  & 515.7439 \\
\rowcolor[HTML]{EFEFEF} 
bic+Frost+SK &
  21.0100 &
  0.6183 &
  \textbf{0.3146} &
  \textbf{0.2528} &
  \textbf{0.4252} &
  \textbf{0.4281} &
  0.4413 &
  \textbf{0.4453} &
  \textbf{454.3244} &
  \textbf{565.3653} \\ \hline
Lee            & \textbf{24.5162} & 0.5698          & 0.4899 & 0.3871 & 0.6622 & 0.6556 & 0.6618          & 0.6555          & 81.5083   & 107.9425 \\
\rowcolor[HTML]{EFEFEF} 
SK+Lee+SK      & 21.2372          & 0.6291          & 0.3282 & 0.2656 & 0.4436 & 0.4498 & \textbf{0.4446} & 0.4492          & 402.3800  & 490.1480 \\
bil+Lee+bil    & 22.6270          & \textbf{0.6622} & 0.3366 & 0.2632 & 0.4549 & 0.4459 & 0.4702          & 0.4625          & 348.7550  & 482.1399 \\
\rowcolor[HTML]{EFEFEF} 
bil+Lee+SK     & 21.9839          & 0.6440          & 0.3228 & 0.2568 & 0.4363 & 0.4349 & 0.4504          & 0.4509          & 412.9327  & 533.0647 \\
bic+Lee+bic    & 22.3082          & 0.6592          & 0.3316 & 0.2599 & 0.4482 & 0.4402 & 0.4656          & 0.4581          & 367.6261  & 505.5940 \\
\rowcolor[HTML]{EFEFEF} 
bic+Lee+SK &
  21.7294 &
  0.6394 &
  \textbf{0.3175} &
  \textbf{0.2540} &
  \textbf{0.4291} &
  \textbf{0.4302} &
  0.4454 &
  \textbf{0.4475} &
  \textbf{438.0506} &
  \textbf{554.5347} \\ \hline
NLM            & \textbf{26.6022} & 0.6708          & 0.2714 & 0.2789 & 0.3670 & 0.4724 & 0.3654          & 0.4721          & 876.5564  & 401.2516 \\
\rowcolor[HTML]{EFEFEF} 
SK+NLM+SK      & 22.4064          & 0.6853          & 0.1942 & 0.2295 & 0.2624 & 0.3887 & 0.2629          & 0.3876          & 3286.2159 & 886.1761 \\
bil+NLM+bil    & 23.3926          & \textbf{0.7363} & 0.1796 & 0.2233 & 0.2428 & 0.3782 & 0.2496          & 0.3902          & 4331.7213 & 937.2723 \\
\rowcolor[HTML]{EFEFEF} 
bil+NLM+SK &
  23.1245 &
  0.6985 &
  \textbf{0.1675} &
  \textbf{0.2213} &
  \textbf{0.2264} &
  \textbf{0.3748} &
  \textbf{0.2326} &
  \textbf{0.3866} &
  \textbf{5731.4312} &
  \textbf{972.4528} \\
bic+NLM+bic    & 23.1992          & 0.7352          & 0.2079 & 0.2252 & 0.2810 & 0.3815 & 0.2900          & 0.3948          & 2400.1149 & 902.5313 \\
\rowcolor[HTML]{EFEFEF} 
bic+NLM+SK     & 22.9176          & 0.6983          & 0.1950 & 0.2232 & 0.2636 & 0.3781 & 0.2719          & 0.3913          & 3103.7562 & 935.0230 \\ \hline
\end{tabular}%
}
\caption{Performance comparison of despeckling in \textit{Cameraman}.}
\label{tab_cameraman}
\end{table}
\begin{table}[tbph]
\centering
\resizebox{\columnwidth}{!}{%
\begin{tabular}{lcccccccccc}
\hline
\multicolumn{1}{c}{} &
   &
   &
  \multicolumn{2}{c}{SI} &
  \multicolumn{2}{c}{SSI} &
  \multicolumn{2}{c}{SMPI} &
  \multicolumn{2}{c}{ENL} \\ \cline{4-11} 
\multicolumn{1}{c}{\multirow{-2}{*}{}} &
  \multirow{-2}{*}{PSNR} &
  \multirow{-2}{*}{SSIM} &
  ROI1 &
  ROI2 &
  ROI1 &
  ROI2 &
  ROI1 &
  ROI2 &
  ROI1 &
  ROI2 \\ \hline
mean &
  \textbf{25.7872} &
  0.5913 &
  0.3811 &
  0.4523 &
  0.6540 &
  0.6894 &
  0.6545 &
  0.6899 &
  100.0283 &
  73.9220 \\
\rowcolor[HTML]{EFEFEF} 
SK+mean+SK &
  23.9522 &
  0.6522 &
  0.3215 &
  0.3623 &
  0.5518 &
  0.5522 &
  \textbf{0.5514} &
  \textbf{0.5534} &
  198.4374 &
  179.2793 \\
bil+mean+bil &
  24.8824 &
  0.6715 &
  0.3185 &
  0.3674 &
  0.5465 &
  0.5600 &
  0.5702 &
  0.5808 &
  195.3104 &
  162.2715 \\
\rowcolor[HTML]{EFEFEF} 
bil+mean+SK &
  24.4559 &
  0.6631 &
  0.3145 &
  0.3602 &
  0.5397 &
  0.5490 &
  0.5630 &
  0.5688 &
  205.3877 &
  175.9014 \\
bic+mean+bic &
  24.8142 &
  \textbf{0.6764} &
  0.3170 &
  0.3650 &
  0.5440 &
  0.5564 &
  0.5696 &
  0.5794 &
  198.0500 &
  165.6548 \\
\rowcolor[HTML]{EFEFEF} 
bic+mean+SK  & 24.3601 & 0.6639 & \textbf{0.3135} & \textbf{0.3578} & \textbf{0.5379} & \textbf{0.5454} & 0.5633 & 0.5675          & \textbf{207.1675} & \textbf{179.5517} \\ \hline
median &
  22.8247 &
  0.4671 &
  0.4454 &
  0.5065 &
  0.7643 &
  0.7720 &
  0.7633 &
  0.7707 &
  54.1133 &
  47.4229 \\
\rowcolor[HTML]{EFEFEF} 
SK+median+SK &
  23.9179 &
  0.6511 &
  0.3297 &
  0.3669 &
  0.5658 &
  0.5592 &
  \textbf{0.5654} &
  \textbf{0.5615} &
  \textbf{179.5337} &
  \textbf{170.0665} \\
bil+median+bil &
  25.0510 &
  0.6641 &
  0.3353 &
  0.3780 &
  0.5754 &
  0.5761 &
  0.6043 &
  0.6011 &
  157.6999 &
  143.7421 \\
\rowcolor[HTML]{EFEFEF} 
bil+median+SK &
  24.5331 &
  0.6643 &
  0.3298 &
  0.3697 &
  0.5660 &
  0.5635 &
  0.5941 &
  0.5870 &
  168.5829 &
  157.3995 \\
bic+median+bic &
  \textbf{25.1019} &
  \textbf{0.6771} &
  0.3307 &
  0.3738 &
  0.5675 &
  0.5697 &
  0.5964 &
  0.5970 &
  166.5098 &
  149.4663 \\
\rowcolor[HTML]{EFEFEF} 
bic+median+SK &
  24.5400 &
  0.6706 &
  \textbf{0.3255} &
  \textbf{0.3655} &
  \textbf{0.5586} &
  \textbf{0.5570} &
  0.5870 &
  0.5830 &
  177.3728 &
  163.8292 \\ \hline
Frost &
  \textbf{25.5160} &
  0.5876 &
  0.3811 &
  0.4524 &
  0.6540 &
  0.6896 &
  0.6545 &
  0.6901 &
  100.0349 &
  73.8472 \\
\rowcolor[HTML]{EFEFEF} 
SK+Frost+SK &
  23.3892 &
  0.6497 &
  0.3216 &
  0.3623 &
  0.5518 &
  0.5522 &
  \textbf{0.5514} &
  \textbf{0.5535} &
  198.4002 &
  179.1907 \\
bil+Frost+bil &
  24.5838 &
  0.6595 &
  0.3185 &
  0.3674 &
  0.5466 &
  0.5600 &
  0.5702 &
  0.5808 &
  195.2250 &
  162.2008 \\
\rowcolor[HTML]{EFEFEF} 
bil+Frost+SK &
  23.9240 &
  0.6597 &
  0.3145 &
  0.3601 &
  0.5397 &
  0.5489 &
  0.5631 &
  0.5686 &
  205.3131 &
  176.0501 \\
bic+Frost+bic &
  24.5156 &
  \textbf{0.6635} &
  0.3171 &
  0.3651 &
  0.5441 &
  0.5565 &
  0.5696 &
  0.5796 &
  198.0207 &
  165.4331 \\
\rowcolor[HTML]{EFEFEF} 
bic+Frost+SK & 23.8353 & 0.6602 & \textbf{0.3135} & \textbf{0.3579} & \textbf{0.5380} & \textbf{0.5456} & 0.5633 & 0.5676          & \textbf{207.0844} & \textbf{179.3481} \\ \hline
Lee &
  24.4192 &
  0.5369 &
  0.4178 &
  0.4771 &
  0.7169 &
  0.7273 &
  0.7175 &
  0.7279 &
  69.2720 &
  59.6753 \\
\rowcolor[HTML]{EFEFEF} 
SK+Lee+SK &
  23.9687 &
  0.6658 &
  0.3234 &
  0.3652 &
  0.5549 &
  0.5567 &
  \textbf{0.5545} &
  \textbf{0.5581} &
  194.0353 &
  173.4976 \\
bil+Lee+bil &
  \textbf{25.6068} &
  0.6951 &
  0.3207 &
  0.3722 &
  0.5503 &
  0.5673 &
  0.5741 &
  0.5886 &
  190.0062 &
  153.9451 \\
\rowcolor[HTML]{EFEFEF} 
bil+Lee+SK &
  24.8063 &
  0.6814 &
  0.3165 &
  0.3640 &
  0.5431 &
  0.5548 &
  0.5665 &
  0.5748 &
  200.3344 &
  168.6571 \\
bic+Lee+bic &
  25.4260 &
  \textbf{0.6960} &
  0.3186 &
  0.3682 &
  0.5467 &
  0.5612 &
  0.5723 &
  0.5846 &
  194.2558 &
  159.9936 \\
\rowcolor[HTML]{EFEFEF} 
bic+Lee+SK   & 24.6456 & 0.6801 & \textbf{0.3147} & \textbf{0.3603} & \textbf{0.5400} & \textbf{0.5491} & 0.5653 & 0.5716          & \textbf{204.0991} & \textbf{174.6738} \\ \hline
NLM &
  \textbf{26.8911} &
  0.6733 &
  0.3190 &
  0.3106 &
  0.5474 &
  0.4735 &
  0.5472 &
  0.4727 &
  204.4260 &
  335.1162 \\
\rowcolor[HTML]{EFEFEF} 
SK+NLM+SK &
  25.0377 &
  0.7171 &
  0.3017 &
  0.3056 &
  0.5178 &
  0.4658 &
  \textbf{0.5177} &
  0.4684 &
  255.2428 &
  352.5519 \\
bil+NLM+bil &
  26.5145 &
  0.7529 &
  0.2936 &
  0.2857 &
  0.5038 &
  0.4355 &
  0.5261 &
  0.4536 &
  270.1142 &
  441.1793 \\
\rowcolor[HTML]{EFEFEF} 
bil+NLM+SK   & 25.7139 & 0.7240 & \textbf{0.2921} & \textbf{0.2794} & \textbf{0.5012} & \textbf{0.4258} & 0.5234 & \textbf{0.4432} & \textbf{275.7648} & \textbf{483.0673} \\
bic+NLM+bic &
  26.6767 &
  \textbf{0.7639} &
  0.3039 &
  0.3270 &
  0.5215 &
  0.4984 &
  0.5472 &
  0.5213 &
  233.9957 &
  255.7811 \\
\rowcolor[HTML]{EFEFEF} 
bic+NLM+SK &
  25.8931 &
  0.7363 &
  0.3016 &
  0.3181 &
  0.5175 &
  0.4849 &
  0.5430 &
  0.5066 &
  241.2933 &
  285.9154 \\ \hline
\end{tabular}%
}
\caption{Performance comparison of despeckling in \textit{Boat}.}
\label{tab_boat}
\end{table}

We also examined an alternative Up-Down scaling approach. Here, the noisy image was first upscaled to $2N\times 2M$, filtered, and then downscaled back to its original dimensions. The results, reported in Table \ref{tab_UpDown} for \textit{Cameraman} using the NLM filter, show a deterioration in all quality metrics. This suggests that upscaling a noisy image before filtering does not improve denoising; instead, it amplifies noise variance, making filtering less effective. Although higher resolution might be expected to improve noise reduction, in practice, it also enlarges the noise pattern, leading to worse results. This confirms that Up-Down scaling is not a viable strategy for speckle noise reduction.\\
\begin{table}[tbph]
\centering
\resizebox{\columnwidth}{!}{%
\begin{tabular}{lcccccc}
\hline
                                              & PSNR    & SSIM   & SI     & SSI    & SMPI   & ENL     \\ \hline
noisy image                                   & 18.6346 & 0.4171 & 0.7398 & -      & -      & 15.6282 \\ \hline
\rowcolor[HTML]{EFEFEF} 
NLM & \textbf{26.5998} & \textbf{0.6705} & \textbf{0.2719} & \textbf{0.3675} & \textbf{0.3659} & \textbf{870.4510} \\ \hline
upscaling with SK + NLM                      & -       & -      & 0.5276 & 1.0000 & 1.0000 & 60.8100 \\
\rowcolor[HTML]{EFEFEF} 
upscaling + NLM + downscaling with SK       & 23.1853 & 0.5544 & 0.5045 & 0.6819 & 0.6811 & 72.5908 \\ \hline
upscaling with bilinear + NLM                & -       & -      & 0.5925 & 1.0000 & 1.0000 & 36.5982 \\
\rowcolor[HTML]{EFEFEF} 
upscaling + NLM + downscaling with bilinear & 22.4021 & 0.5087 & 0.5758 & 0.7783 & 0.8048 & 40.6834 \\ \hline
upscaling with bicubic + NLM                 & -       & -      & 0.5920 & 1.0000 & 1.0000 & 36.7104 \\
\rowcolor[HTML]{EFEFEF} 
upscaling + NLM + downscaling with bicubic  & 22.8475 & 0.5243 & 0.5538 & 0.7485 & 0.7758 & 47.4158 \\ \hline
\end{tabular}%
}
\caption{Results of Up-Down scaling method for \textit{Cameraman} (ROI1$=[220,200,30,40]$) using NLM filter.}
\label{tab_UpDown}
\end{table}

\begin{figure}[htbp]
    \centering
    \begin{subfigure}{0.22\textwidth}
        \includegraphics[width=\linewidth]{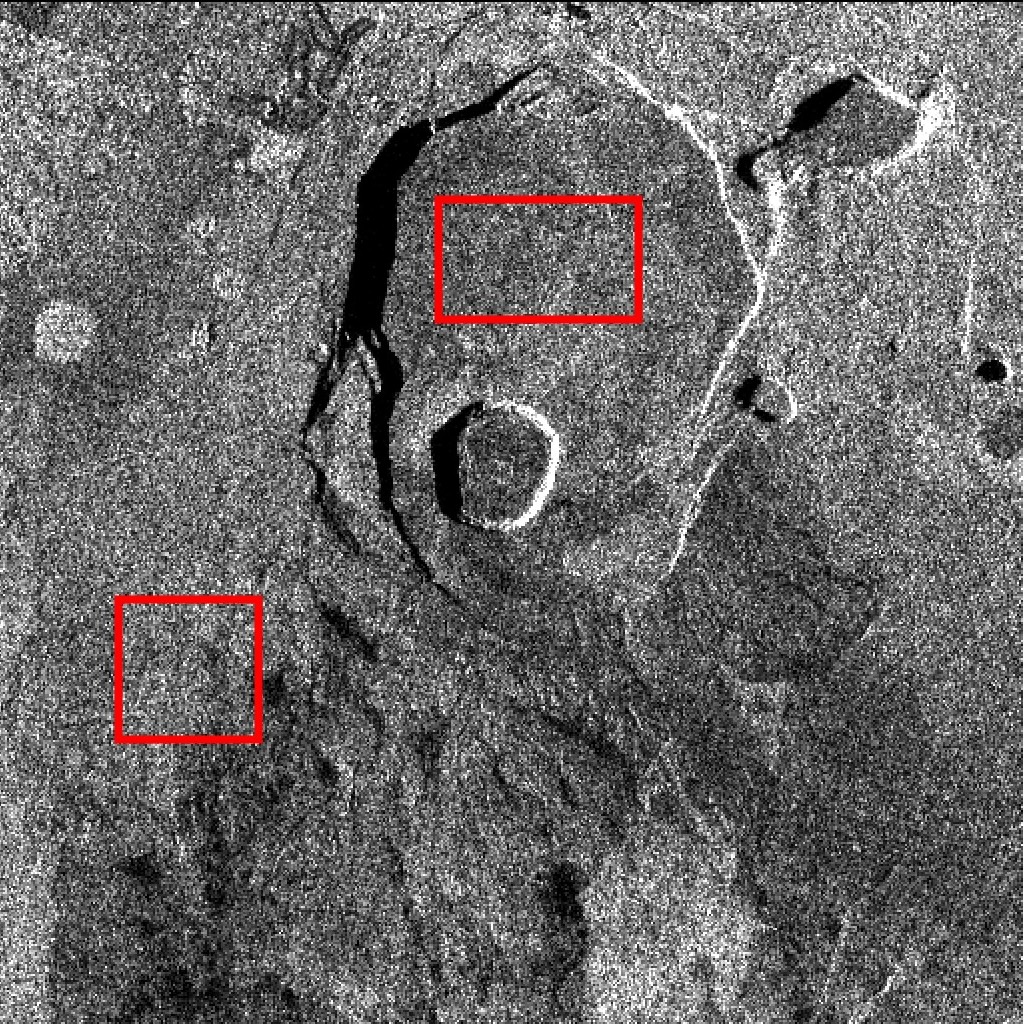}
        \caption{}
    \end{subfigure}
    \begin{subfigure}{0.22\textwidth}
        \includegraphics[width=\linewidth]{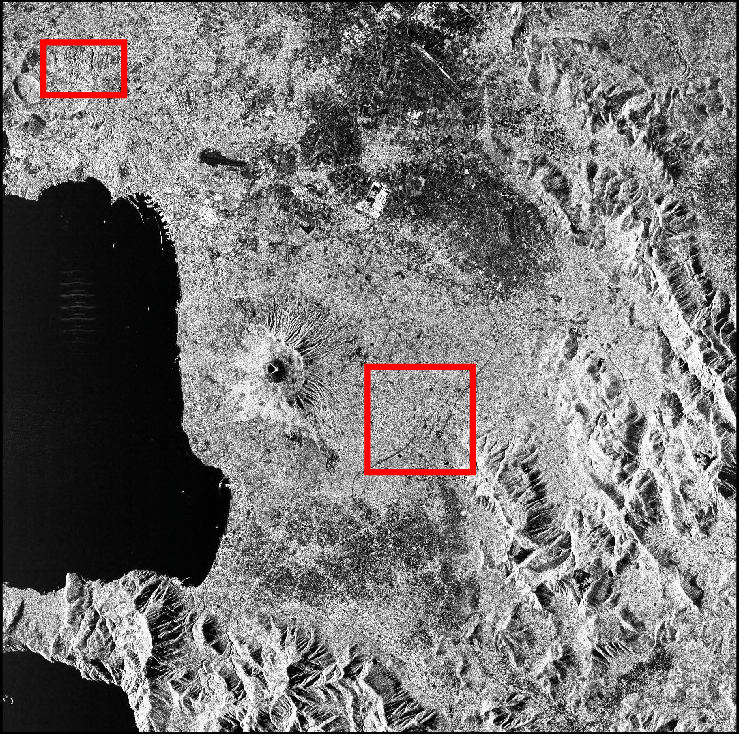}
        \caption{}
    \end{subfigure}
    \begin{subfigure}{0.22\textwidth}
        \includegraphics[width=\linewidth]{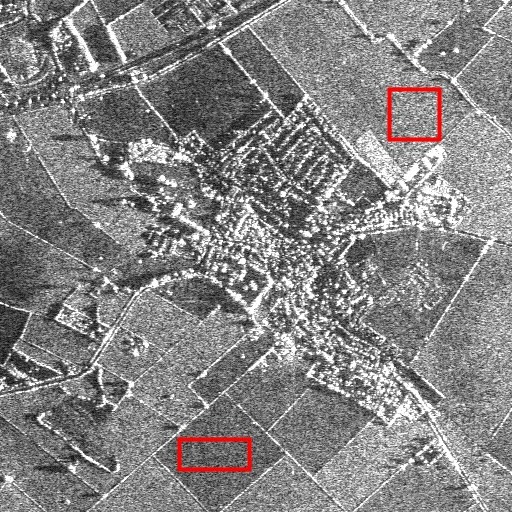}
        \caption{}  
    \end{subfigure}
\caption{(a) \textit{SAR image 1} \cite{url1}: ROI1$=[120,600,140,140]$, ROI2$=[440,200,100,120]$ ; (b) \textit{SAR image 2} \cite{url4}: ROI1$=[900, 900 , 260, 260]$, ROI2$=[100, 100, 200, 130]$; (c) \textit{SAR image 3} \cite{url2}: ROI1$=[390, 90, 50, 50]$, ROI2$=[180, 440, 70, 30]$.}
\label{moon}
\end{figure}
Finally, we test three real SAR images (see Figure \ref{moon}). The real SAR image 1 captures the Kilauea volcano in Hawaii  [$256\times 256$] \cite{url1}; the real SAR image 2 depicts the Vesuvius volcano in Italy [$1800\times 1800$] \cite{url4}. Both were acquired using the Spaceborne Imaging Radar-C/X-band Synthetic Aperture Radar. Finally, the real SAR image 3 depicts a rural scene in Bedfordshire [$512\times 512$] \cite{new, url2}. No prior speckle filtering was applied to these real SAR images, they show the characteristic granular noise inherent to SAR amplitude data. We evaluated the filtering methods both on their own and in combination with the proposed Down-Up processing approach. The results for direct filtering and the best-performing configuration, which combines bicubic downscaling, filtering, and SK upscaling, are presented in Tables \ref{tab_real1}-\ref{tab_real3}. Since the noise distribution in the the real SAR data is not known, the real effectiveness and
strength the despeckling methods can not be reliably measured on them alone. For this reason, we also tested the algorithms on optical test images (\textit{Cameraman} and \textit{Boat}) with artificially added speckle noise. Although these images differ from real SAR data -- since SAR scenes often contain large open areas and strong clutter caused by intense backscattering -- the experiments allow a more reliable evaluation of despeckling performance under known noise conditions. As with the previous tests, the proposed approach consistently improves all quality metrics, highlighting its practical benefits for SAR image denoising.

\begin{table}[tbph]
\centering
\resizebox{\columnwidth}{!}{%
\begin{tabular}{lcccccccc}
\hline
& \multicolumn{2}{c}{SI} & \multicolumn{2}{c}{SSI} & \multicolumn{2}{c}{SMPI} & \multicolumn{2}{c}{ENL} \\ \cline{2-9} 
              & ROI1   & ROI2   & ROI1   & ROI2   & ROI1   & ROI2   & ROI1    & ROI2    \\ \hline
mean          & 0.8551 & 0.8378 & 0.8159 & 0.7971 & 0.8159 & 0.7971 & 10.5791 & 12.1950 \\
\rowcolor[HTML]{EFEFEF} 
bic+mean+SK   & \textbf{0.6669} & \textbf{0.6296} & \textbf{0.6363} & \textbf{0.5990} & \textbf{0.6669} & \textbf{0.6271} & \textbf{26.7768} & \textbf{35.8399} \\ \hline
median        & 0.9252 & 0.9133 & 0.8828 & 0.8689 & 0.8730 & 0.8588 & 8.0196  & 8.9825  \\
\rowcolor[HTML]{EFEFEF} 
bic+median+SK & \textbf{0.6990} & \textbf{0.6676} & \textbf{0.6670} & \textbf{0.6351} & \textbf{0.6654} & \textbf{0.6331} & \textbf{23.8845} & \textbf{30.5810} \\ \hline
Frost         & 0.8635 & 0.8472 & 0.8239 & 0.8060 & 0.8229 & 0.8049 & 10.2192 & 11.7172 \\
\rowcolor[HTML]{EFEFEF} 
bic+Frost+SK  & \textbf{0.6685} & \textbf{0.6314} & \textbf{0.6379} & \textbf{0.6007} & \textbf{0.6684} & \textbf{0.6287} & \textbf{26.5189} & \textbf{35.4536} \\ \hline
Lee           & 0.9269 & 0.9178 & 0.8844 & 0.8731 & 0.8861 & 0.8748 & 7.6431  & 8.4476  \\
\rowcolor[HTML]{EFEFEF} 
bic+Lee+SK    & \textbf{0.6953} & \textbf{0.6614} & \textbf{0.6634} & \textbf{0.6293} & \textbf{0.6962} & \textbf{0.6597} & \textbf{22.6169} & \textbf{29.3665} \\ \hline
NLM           & 1.0439 & 1.0452 & 0.9961 & 0.9944 & 0.9961 & 0.9943 & 4.7630  & 5.0366  \\
\rowcolor[HTML]{EFEFEF} 
bic+NLM+SK    & \textbf{0.6042} & \textbf{0.5292} & \textbf{0.5765} & \textbf{0.5035} & \textbf{0.5891} & \textbf{0.5148} & \textbf{41.1697} & \textbf{74.2178} \\ \hline
\end{tabular}%
}
\caption{Performance comparison of despeckling in \textit{SAR image 1}.}
\label{tab_real1}
\end{table}
\begin{table}[tbph]
\centering
\resizebox{\columnwidth}{!}{%
\begin{tabular}{lcccccccc}
\hline
 & \multicolumn{2}{c}{SI} & \multicolumn{2}{c}{SSI} & \multicolumn{2}{c}{SMPI} & \multicolumn{2}{c}{ENL} \\ \cline{2-9} 
              & ROI1   & ROI2   & ROI1   & ROI2   & ROI1   & ROI2   & ROI1     & ROI2      \\ \hline
mean          & 0.5902 & 0.6423  & 0.7937 & 0.8133 & 0.7937 & 0.8133 & 21.8348 & 15.4842 \\ 
\rowcolor[HTML]{EFEFEF} 
bic+mean+SK   & \textbf{0.4595} & \textbf{0.4901} & \textbf{0.6180} & \textbf{0.6206} & \textbf{0.6578} & \textbf{0.6675} & \textbf{54.9918} & \textbf{41.7461} \\ \hline
median        & 0.6289 & 0.6902 & 0.8457 & 0.8740 & 0.8499 &  0.8862 & 16.8360 & 11.4143 \\ 
\rowcolor[HTML]{EFEFEF} 
bic+median+SK & \textbf{0.4663} & \textbf{0.4997} & \textbf{0.6271} & \textbf{0.6328} & \textbf{0.6734} & \textbf{0.6913} & \textbf{51.3041} & \textbf{37.8912}  \\ \hline
Frost         & 0.5930 & 0.6474 & 0.7975 & 0.8197 & 0.7973 & 0.8193 & 21.4495 & 15.0453 \\ 
\rowcolor[HTML]{EFEFEF} 
bic+Frost+SK  & \textbf{0.4605} & \textbf{0.4912} & \textbf{0.6192} & \textbf{0.6220} & \textbf{0.6590} & \textbf{0.6688} & \textbf{54.5725} & \textbf{41.3903} \\ \hline 
Lee  & 0.6159 & 0.6751 & 0.8283 & 0.8548 & 0.8282 & 0.8547 & 18.4128 & 12.6953 \\ 
\rowcolor[HTML]{EFEFEF} 
bic+Lee+SK    & \textbf{0.4666} & \textbf{0.5003} & \textbf{0.6274} & \textbf{0.6335} & \textbf{0.6675} & \textbf{0.6810} & \textbf{51.7873} & \textbf{38.4733} \\ \hline
NLM           & 0.6643 & 0.7435 & 0.8933 & 0.9415 & 0.8961 &  0.9435 & 13.5516 & 8.5995 \\ 
\rowcolor[HTML]{EFEFEF} 
bic+NLM+SK    & \textbf{0.4760} & \textbf{0.5286} & \textbf{0.6401} & \textbf{0.6693} & \textbf{0.6835} & \textbf{0.7213} & \textbf{47.5864} & \textbf{30.7837} \\ \hline
\end{tabular}%
}
\caption{Performance comparison of despeckling in \textit{SAR image 2}.}
\label{tab_real2}
\end{table}
\begin{table}[tbph]
\centering
\resizebox{\columnwidth}{!}{%
\begin{tabular}{lcccccccc}
\hline
 & \multicolumn{2}{c}{SI} & \multicolumn{2}{c}{SSI} & \multicolumn{2}{c}{SMPI} & \multicolumn{2}{c}{ENL} \\ \cline{2-9} 
              & ROI1   & ROI2   & ROI1   & ROI2   & ROI1   & ROI2   & ROI1     & ROI2      \\ \hline
mean          & 0.4526 & 0.5740 & 0.6998 & 0.6902 & 0.7004 & 0.6903 & 68.6137  & 59.7535   \\
\rowcolor[HTML]{EFEFEF} 
bic+mean+SK   & \textbf{0.3050} & \textbf{0.3796} & \textbf{0.4716} & \textbf{0.4564} & \textbf{0.4914} & \textbf{0.4764} & \textbf{316.3150} & \textbf{293.9065}  \\ \hline
median        & 0.4909 & 0.6331 & 0.7590 & 0.7612 & 0.7566 & 0.7648 & 50.4170  & 40.1180   \\
\rowcolor[HTML]{EFEFEF} 
bic+median+SK & \textbf{0.3166} & \textbf{0.3960} & \textbf{0.4896} & \textbf{0.4761} & \textbf{0.5075} & \textbf{0.4965} & \textbf{274.0890} & \textbf{248.4684}  \\ \hline
Frost         & 0.4533 & 0.5748 & 0.7009 & 0.6911 & 0.7014 & 0.6912 & 68.2022  & 59.4248   \\
\rowcolor[HTML]{EFEFEF} 
bic+Frost+SK  & \textbf{0.3052} & \textbf{0.3799} & \textbf{0.4720} & \textbf{0.4567} & \textbf{0.4917} & \textbf{0.4767} & \textbf{315.3717} & \textbf{293.0303}  \\ \hline
Lee           & 0.4991 & 0.6052 & 0.7716 & 0.7277 & 0.7722 & 0.7277 & 46.4238  & 48.3649   \\
\rowcolor[HTML]{EFEFEF} 
bic+Lee+SK    & \textbf{0.3158} & \textbf{0.3872} & \textbf{0.4883} & \textbf{0.4656} & \textbf{0.5085} & \textbf{0.4859} & \textbf{275.3650} & \textbf{271.4871}  \\ \hline
NLM           & 0.4705 & 0.4136 & 0.7274 & 0.4973 & 0.7271 & 0.4969 & 58.9538  & 222.3438  \\
\rowcolor[HTML]{EFEFEF} 
bic+NLM+SK    & \textbf{0.2390} & \textbf{0.2254} & \textbf{0.3695} & \textbf{0.2711} & \textbf{0.3848} & \textbf{0.2865} & \textbf{839.9393} & \textbf{2320.4277} \\ \hline
\end{tabular}%
}
\caption{Performance comparison of despeckling in \textit{SAR image 3}.}
\label{tab_real3}
\end{table}
\section{Conclusions and final remarks}\label{s5}
In this paper, we proposed two methods for image restoration problems motivated by remote sensing applications, namely missing data and speckle noise. The first method focuses on gap filling using the new introduced LP-SK algorithm, which relies only on past sample values and does not require auxiliary or future data. This method could be suited to remote-sensing applications, where data are often acquired progressively along spatial tracks, as in GNSS-R scanning \cite{Zavo}. Similar scenarios arise in airborne remote sensing, radar imaging, or even medical imaging, where full data acquisition may be delayed or interrupted, and timely reconstruction is crucial, see e.g. \cite{ex1, ex2, ex3}.\\
The second method deals with speckle noise reduction. We introduced a Down-Up scaling approach that includes three steps: downscaling with bicubic interpolation, classical despeckling filtering, and upscaling with SK algorithm. This combination provided good results in terms of noise suppression. The experiments on standard test images provide a controlled validation under known ground-truth conditions, whereas the tests on SAR data show the potential of the proposed approach in a remote sensing setting.\\

The algorithms developed in this paper are part of the objectives of the research project PRIN 2022 named ``AIDA'' (``AI- and DIP-Enhanced DAta Augmentation for Remote Sensing of Soil Moisture and Forest Biomass''), funded by the European Union under the Italian National Recovery and Resilience Plan (NRRP) of NextGenerationEU, under the Italian Ministry of Universities and Research.\\

We also test both methods together on a real SAR image (\textit{SAR image 1}, see Figure \ref{moon} (a)). As shown in Figure \ref{final_test}, we first applied the LP-SK algorithm to fill in missing data, and then applied the Down-Up despeckling procedure using the NLM filter. In Figure \ref{zoom}, a zoomed-in comparison of selected ROIs is provided. The result of this preliminary test is promising and suggests that this combined approach can be effective. In future work, we plan to extend the application of LP-SK algorithm to real satellite images affected by missing data, with the goal of improving the reliability of remote sensing analyzes in practical scenarios.\\
\begin{figure}[htbp]
    \centering
    \begin{subfigure}{0.22\textwidth}
        \includegraphics[width=\linewidth]{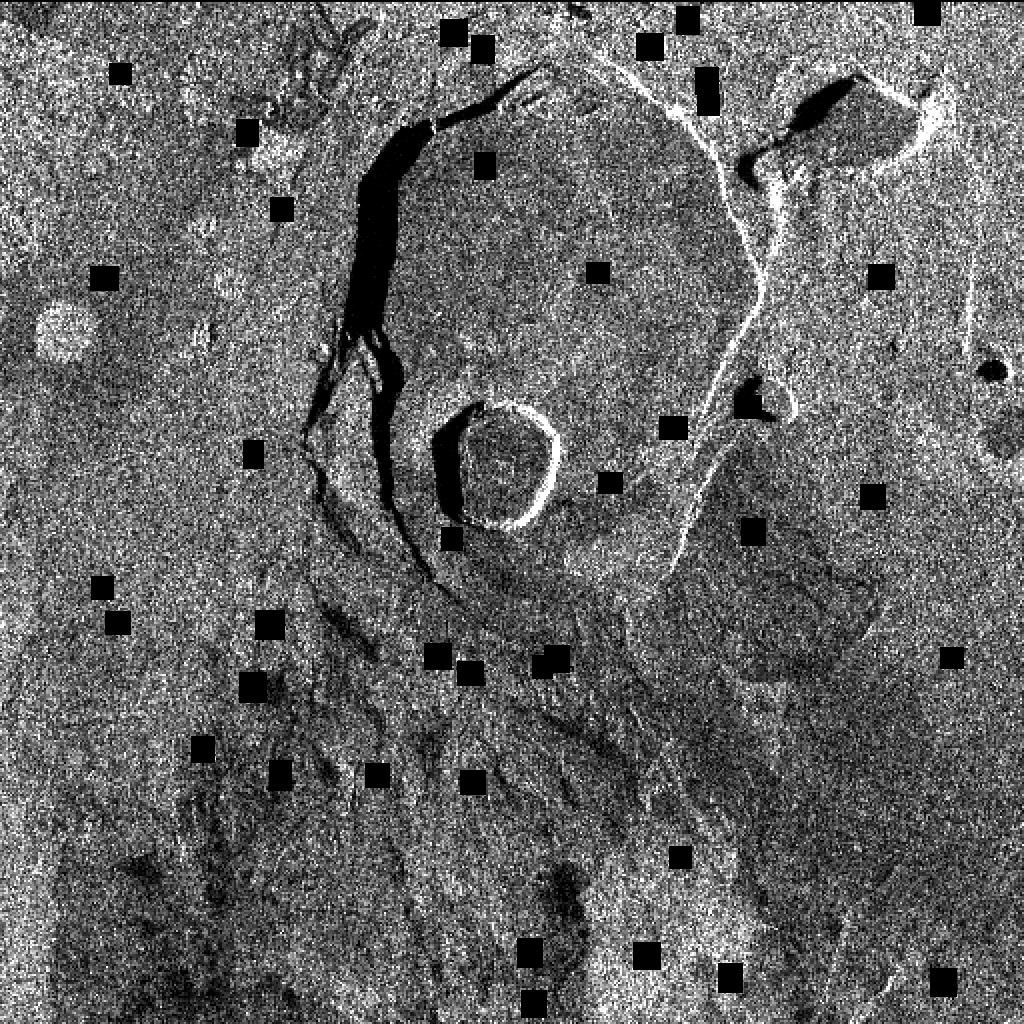}
        \caption{}
    \end{subfigure}
    \begin{subfigure}{0.22\textwidth}
        \includegraphics[width=\linewidth]{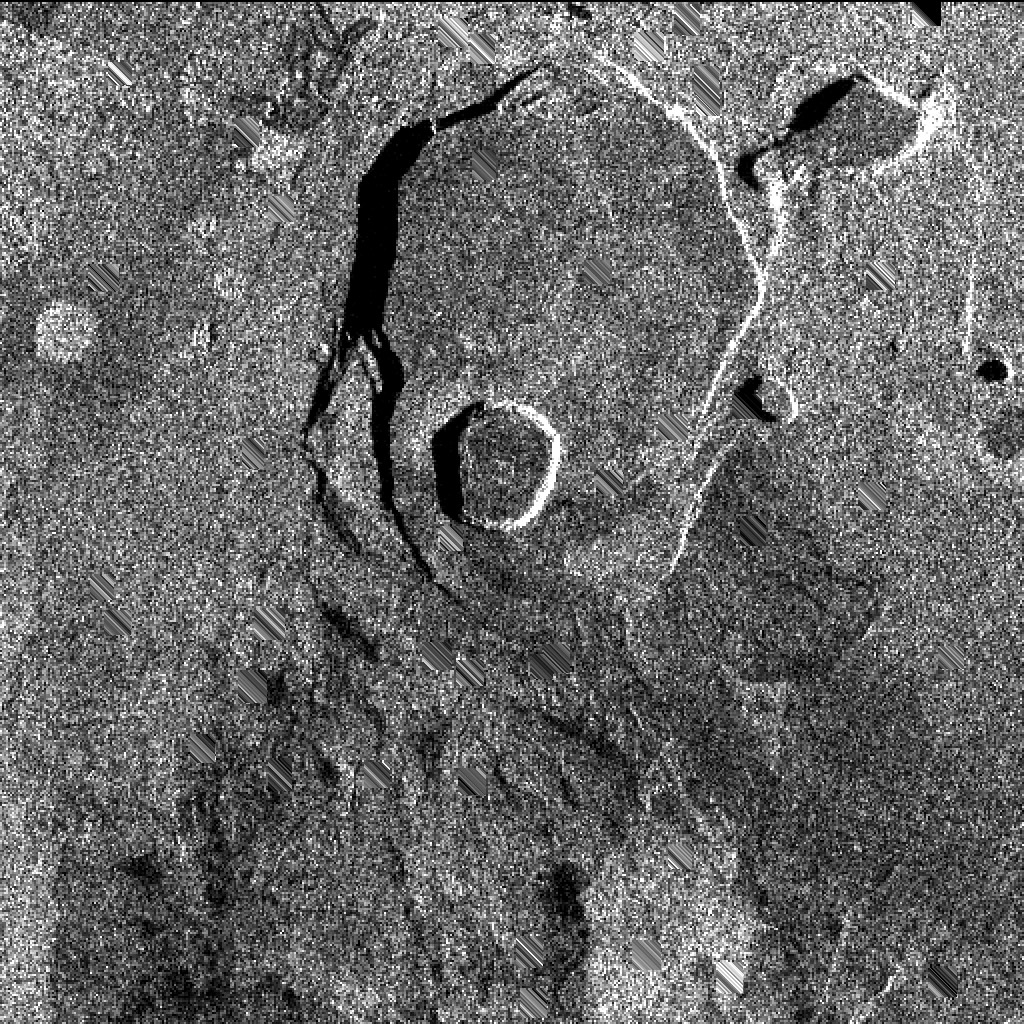}
        \caption{}
    \end{subfigure}
    \begin{subfigure}{0.22\textwidth}
        \includegraphics[width=\linewidth]{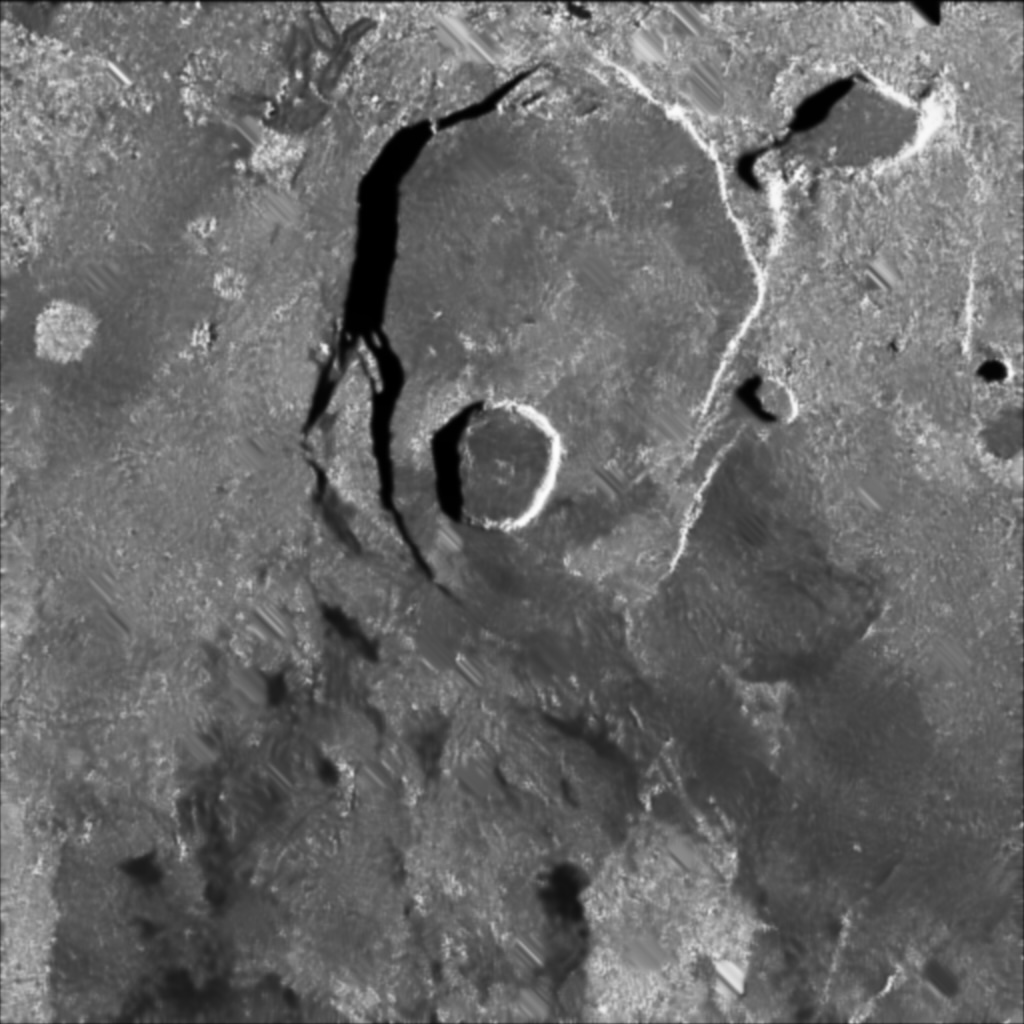}
        \caption{}
    \end{subfigure}
\caption{\textit{SAR image 1}: (a) image with simulated loss ($2.53$\% missing pixels); (b) reconstructed image with SK linear prediction; (c) despeckled image obtained with the proposed Down-Up scaling method.}
\label{final_test}
\end{figure}

\begin{figure}[htbp]
    \centering
    \begin{subfigure}{0.22\textwidth}
        \includegraphics[width=\linewidth]{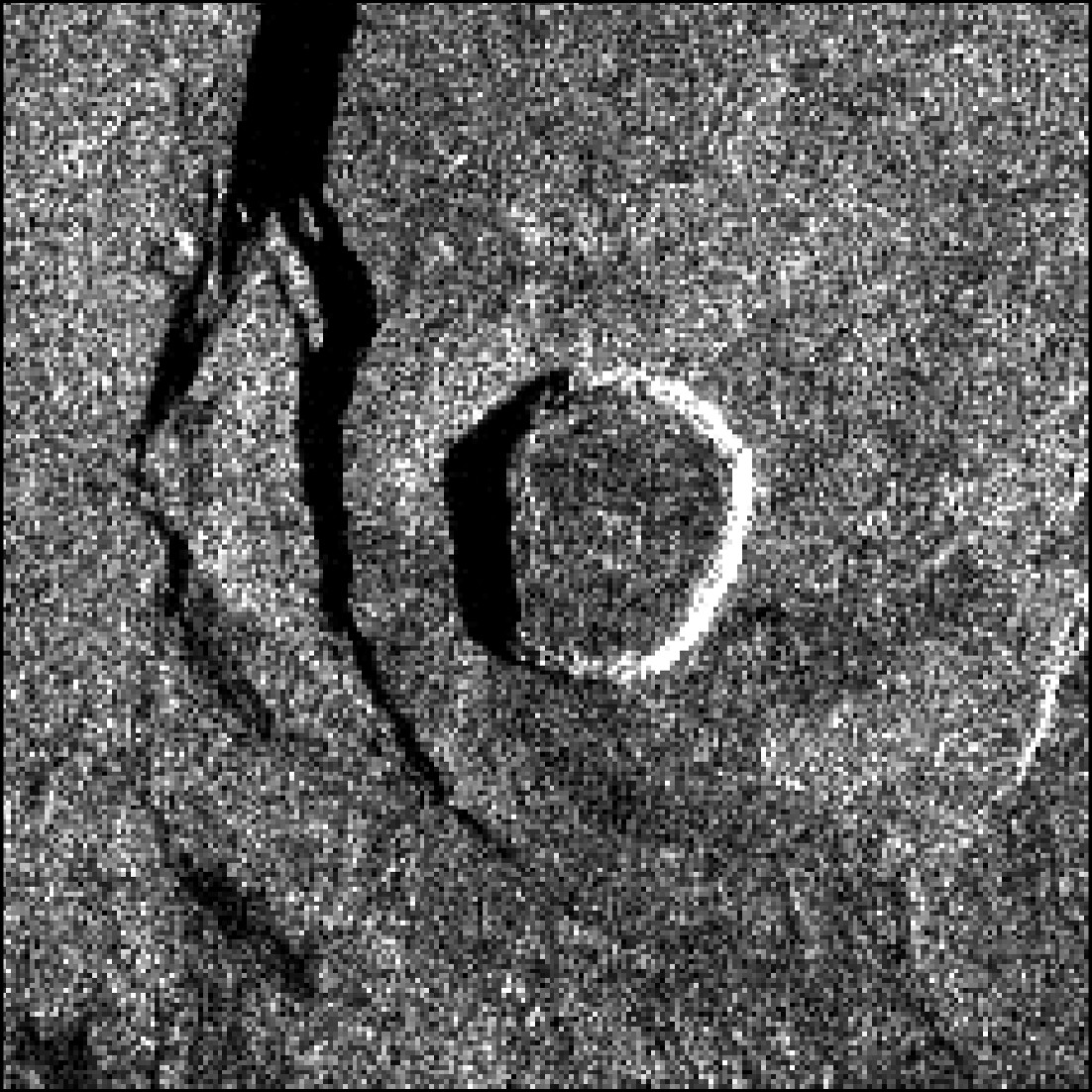}
        \caption{}
    \end{subfigure}
    \begin{subfigure}{0.22\textwidth}
        \includegraphics[width=\linewidth]{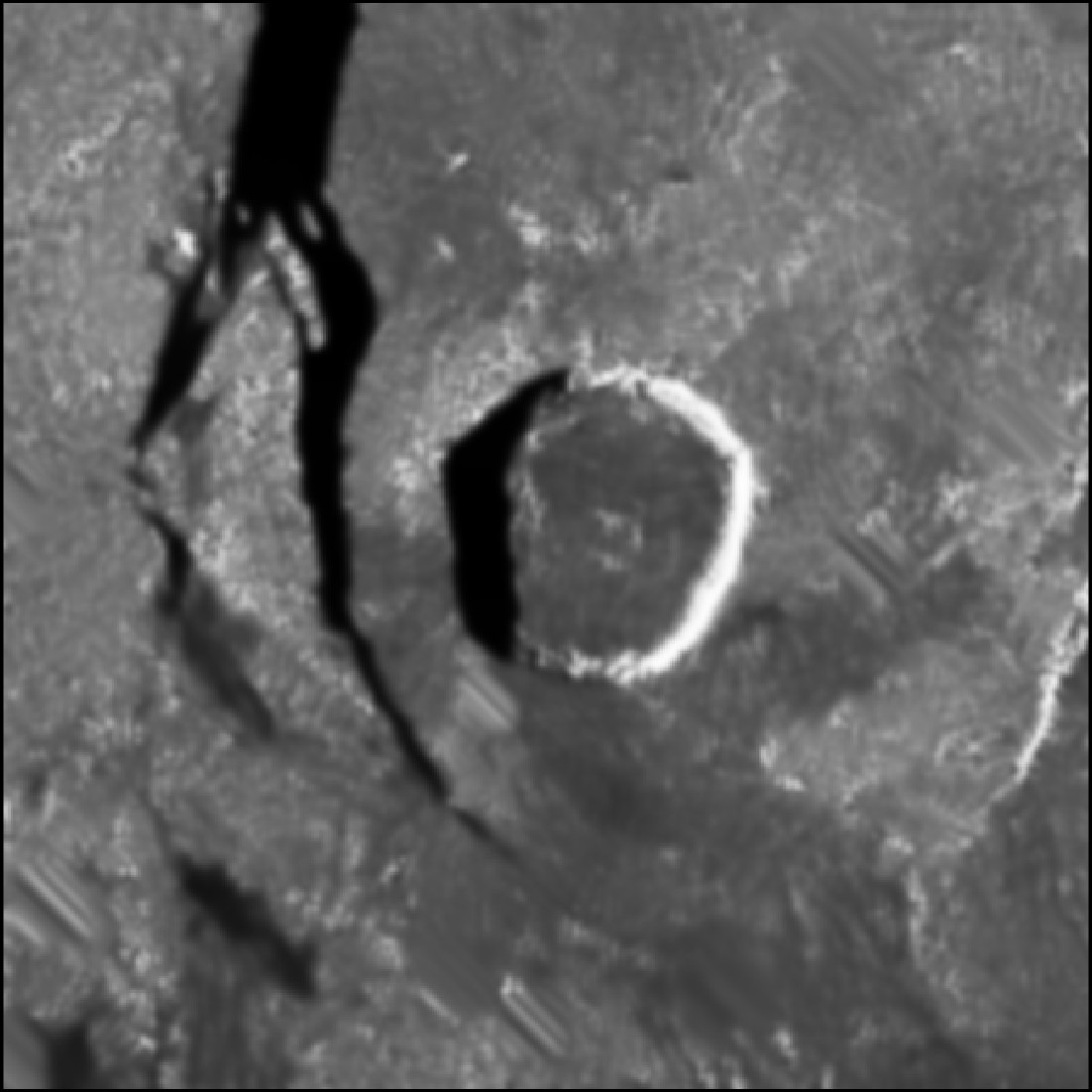}
        \caption{}
    \end{subfigure}
    \begin{subfigure}{0.22\textwidth}
        \includegraphics[width=\linewidth]{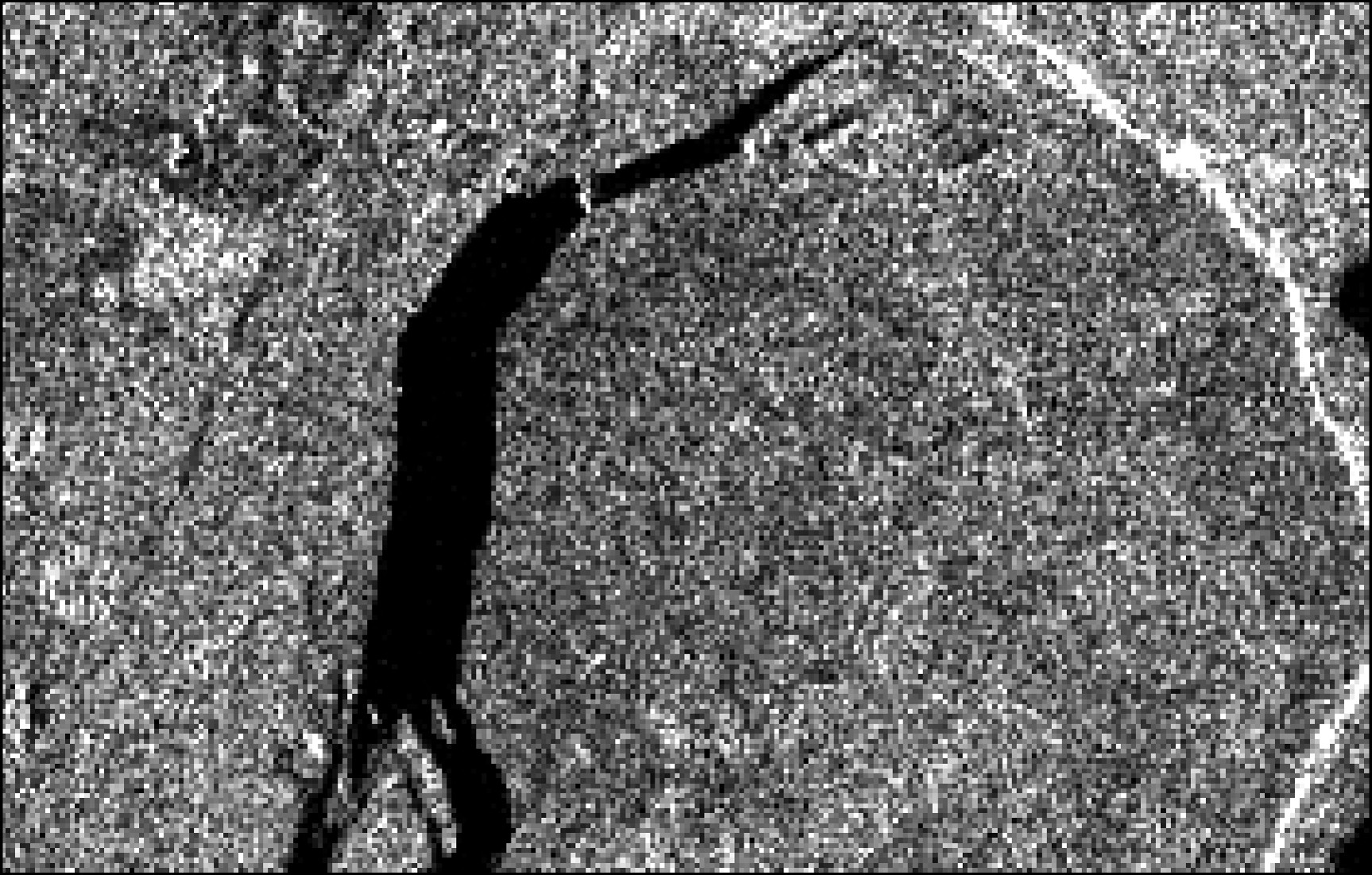}
        \caption{}
    \end{subfigure}
    \begin{subfigure}{0.22\textwidth}
        \includegraphics[width=\linewidth]{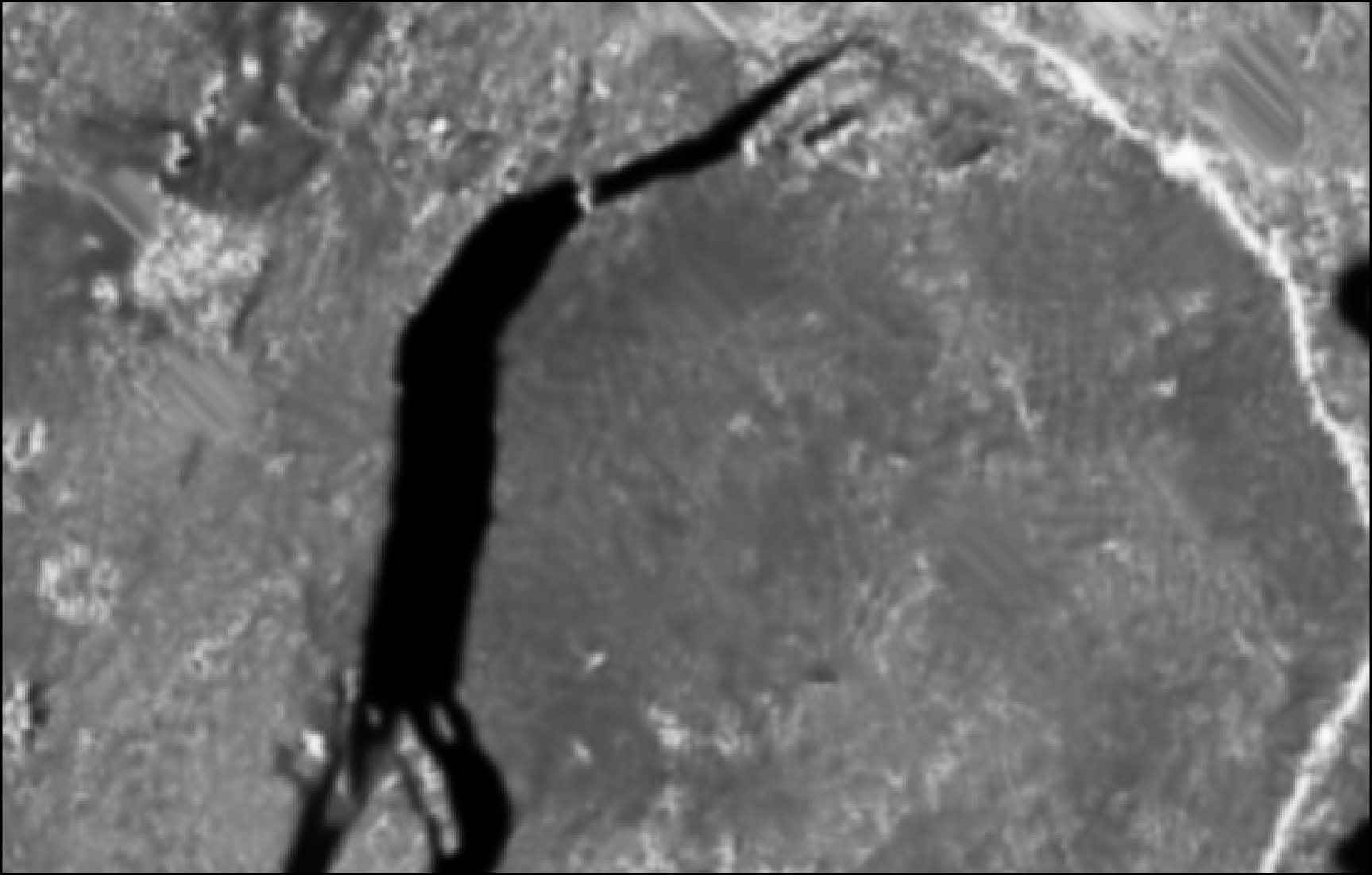}
        \caption{}
    \end{subfigure}
\caption{Zoomed-in comparison on selected ROIs for \textit{SAR image 1}: (a)-(b) central crater region [250, 250, 450, 450]; (c)-(d) fracture on the left side [300, 120, 220, 380].}
\label{zoom}
\end{figure}

Finally, we can also stress that the theoretical results about the order of approximation of SK operators are crucial to analyze the accuracy performances of both the SK and LP-SK algorithms. In particular, Corollary \ref{corollary_Lipschitz} is strictly related to the decays of some of the numerical errors computed in Section \ref{s3} and Section \ref{s4}.

\subsubsection*{Acknowledgments} 
{\small The authors are members of the Gruppo Nazionale per l'Analisi Matematica, la Probabilit\`{a} e le loro Applicazioni (GNAMPA) of the Istituto Nazionale di Alta Matematica (INdAM), of the gruppo UMI (Unione Matematica Italiana) T.A.A. (Teoria dell' Approssimazione e Applicazioni), and of the network RITA (Research ITalian network on Approximation). }

\subsubsection*{Funding}
{\small The authors have been supported within the project PRIN 2022: ``AI- and DIP-Enhanced DAta Augmentation for Remote Sensing of Soil Moisture and Forest Biomass (AIDA)'' funded by the European Union under the Italian National Recovery and Resilience Plan (NRRP) of NextGenerationEU, under the Italian Ministry of Universities and Research (Project Code: 20229FX3B9, CUP: J53D23000660001).}
\subsubsection*{Disclosure of Interests}
{\small The authors declare that they have no conflict of interest and competing interest.}
\subsubsection*{Availability of data and material and Code availability}
{\small All the data generated for this study were stored in our laboratory and are not publicly available. Researchers who wish to access the data directly contacted the corresponding author.}

%
% ---- Bibliography ----
%

\end{document}